\documentclass{article}

\usepackage{a4,latexsym,amssymb,amsmath,amsthm,paralist}
\usepackage[T1]{fontenc}
\usepackage{lmodern}

\usepackage[all]{xy}
\entrymodifiers={+!!<0pt,\fontdimen22\textfont2>}
\def\mapr#1{\stackrel{#1}{\longrightarrow}}

\def\surjd#1{\lower4pt\hbox{$\downarrow$}\kern-5.65pt\Big\downarrow\rlap {$\vcenter{\hbox{$\scriptstyle{{#1}}$}}$}}

\newcommand{\Spec}{\text{\it Spec}}
\newcommand{\Q}{{\mathbb Q}}
\newcommand{\Z}{{\mathbb Z}}
\newcommand{\F}{{\mathbb F}}
\newcommand{\Frob}{{\mathit{Frob}}}
\newcommand{\N}{{\mathbb N}}
\newcommand{\G}{{\mathbb G}}
\renewcommand{\O}{{\mathcal{O}}}

\newcommand{\Cl}{\text{\it Cl}}
\newcommand{\p}{{\mathfrak p}}
\newcommand{\q}{{\mathfrak q}}
\newcommand{\sm}{{\,\smallsetminus\,}}
\newcommand{\et}{\mathit{et}}
\newcommand{\nrel}{\textit{nr\!,el}}
\newcommand{\fl}{\text{\it fl}\,}
\newcommand{\cd}{\text{\it cd}\:}

\newcommand{\Gal}{\mathit{G}}

\newcommand{\Hom}{\text{\rm Hom}}
\font\emas = cmsy10 scaled\magstep2
\font\smallemas = cmsy10
\newcommand{\freeproductmed}{\mathop{\lower.2mm\hbox{\emas \symbol{3}}}\limits}
\newcommand{\freeproductsmall}{\mathop{\lower.2mm\hbox{\smallemas \symbol{3}}}}
\newcommand{\lang}{\longrightarrow}
\font\russ=wncyr10
\renewcommand{\min}{\text{\rm min}}
 \newcommand{\ressum}{\mathop{\hbox{${\displaystyle\bigoplus}'$}}\limits}
 \newcommand{\ressumsmall}{\mathop{\hbox{${\bigoplus}'$}}}
 \newcommand{\coker}{\text{\rm coker}}
\newcommand{\ds}{\displaystyle}
\renewcommand{\a}{\mathfrak{a}}

\newcommand{\nr}{\mathit{nr}}
\newcommand{\el}{\mathit{el}}
\renewcommand{\P}{\mathfrak{P}}
\newcommand{\im}{\mathrm{im}}
\newcommand{\M}{\mathcal{M}}
\newcommand{\Et}{\mathrm{Et}}
\newcommand{\FEt}{\mathrm{FEt}}
\newcommand{\Pic}{\mathrm{Pic}}
\newcommand{\Br}{\mathrm{Br}}
\newcommand{\Mor}{\mathrm{Mor}}
\newcommand{\liso}{\stackrel{\sim}{\lang}}
\newcommand{\rec}{\mathit{rec}}
\newcommand{\T}{\mathcal{T}}
\newcommand{\mQ}{\mathfrak{Q}}
\newcommand{\cor}{\mathit{cor}}

\font\russ=wncyr10
 1

\def\Sha{\hbox{\russ\char88}}

\def\Be{\hbox{\russ\char66}}

\newcommand{\back}{\hspace{-1em}}

\newtheoremstyle{alex}
  {}
  {}
  {\sl}
  {}
  {\bf}
  {.}
  {.5em}
  {}
\newtheoremstyle{alexdef}
  {}
  {}
  {\rm }
  {}
  {\bf}
  {.}
  {.5em}
  {}
\theoremstyle{alex}
 
\newtheorem{theorem}{Theorem}[section]
\newtheorem{corollary}[theorem]{Corollary}
\newtheorem{lemma}[theorem]{Lemma}

\newtheorem{proposition}[theorem]{Proposition}

\theoremstyle{alexdef}
\newtheorem{definition}[theorem]{Definition}
\newtheorem{remark}[theorem]{Remark}
\newtheorem{example}[theorem]{Example}
\newtheorem{remarks}[theorem]{Remarks}

\title{\bf\boldmath On pro-$p$ fundamental groups of \\ marked  arithmetic curves}
\author{Alexander Schmidt}
\date{January 16, 2009}

\begin{document}
\hyphenation{de-fi-ni-tion}
\maketitle
\begin{abstract}
Let $k$ be a global field, $p$ an odd prime number different from $\text{char}(k)$ and $S$, $T$ disjoint, finite sets of primes of $k$. Let $G_S^T(k)(p)=\Gal(k_S^T(p)|k)$ be the Galois group of the maximal $p$-extension of $k$ which is unramified outside $S$ and completely split at $T$.  We prove the existence of a finite set of primes $S_0$, which can be chosen disjoint from any given set $\M$ of Dirichlet density zero, such that the cohomology of $G_{S\cup S_0}^T(k)(p)$ coincides with the \'{e}tale cohomology of the associated marked arithmetic curve. In particular, \hbox{$\cd\, G_{S\cup S_0}^T(k)(p)=2$}. Furthermore, we can choose $S_0$ in such a way that $k_{S\cup S_0}^T(p)$ realizes the maximal $p$-extension $k_\p(p)$ of the local field $k_\p$ for all $\p\in S\cup S_0$, the cup-product
$
H^1(G_{S\cup S_0}^T(k)(p),\F_p) \otimes H^1(G_{S\cup S_0}^T(k)(p),\F_p) \to H^2(G_{S\cup S_0}^T(k)(p),\F_p)
$
is surjective and the decomposition groups of the primes in $S$ establish a free product inside $G_{S\cup S_0}^T(k)(p)$. This generalizes previous work of the author where similar results were shown in the case $T=\varnothing$ under the restrictive assumption  $p\nmid \# \Cl(k)$ and $\zeta_p\notin k$.
\end{abstract}

\section{Introduction}

Let $k$ be a number field,  $S$ a finite set of places of $k$ and $p$ a prime number. We denote the maximal $p$-extension of $k$ which is unramified outside $S$ by $k_S(p)$ and put
\[
G_S(k)(p)=\Gal(k_S(p)|k).
\]
The group $G_S(k)(p)$ is rather well (but far from being completely) understood in the case that $S$ contains all primes dividing $p$ (see Chapter~X of \cite{NSW} for an overview on the results known).
Besides the fact that the group $G_S(k)(p)$ is finitely presentable, until recently not much was known in the case that $S$ does not contain all places above $p$.

In 2005, J.~Labute \cite{La}  considered the case $k=\Q$ and found the first examples of pairs $(S,p)$ such that  $G_S(\Q)(p)$  has cohomological dimension~$2$ and $p$ is not in~$S$. This progress was possible by using the newly developed theory of mild pro-$p$-groups.   The author then generalized  Labute's results to arbitrary number fields \cite{circular,kpi1}, whereas the focus was on the $K(\pi,1)$-property.  This property says that the cohomology of the group $G_S(k)(p)$ coincides with the \'{e}tale cohomology of the scheme $\Spec(\O_k)\sm S$, in particular, it implies cohomological dimension~$2$.  However, for a given number field $k$, finitely many prime numbers had to be excluded:  the divisors of the class number and such $p$ with  $\zeta_p\in k$.
It has turned out now that this restriction can be removed if one considers a more general problem from the very beginning, namely the question for the group  $G_S^T(k)(p)$. Here $S$ and $T$ are finite sets of places,
$k_S^T(p)$ is the maximal $p$-extension of $k$ which is unramified outside $S$ and completely split at~$T$, and $G_S^T(k)(p)=\Gal(k_S^T(p)|k)$. This group is also interesting in the function field case.

\medskip
In this paper we prove without assumptions on $S$ and $T$ that, by adding finitely many primes to $S$, we achieve a situation where $G_S^T(k)(p)$ has cohomological dimension~$2$  and further nice properties. Moreover, we can avoid an arbitrarily given set of Dirichlet density zero (in particular, the places above $p$) when choosing the additional primes. The precise result is the following:

\begin{theorem} \label{haupt} Let $k$ be a global field and let $p$ be an odd prime number different from $\mathrm{char}(k)$.  Let  $S$, $T$ and $\M$ be pairwise disjoint sets of places of~$k$, where $S$ and $T$ are finite and $\M$ has Dirichlet density  $\delta(\M)=0$. Then there exists a finite set of places $S_0$ of $k$ which is disjoint from $S\cup T\cup \M$ and such that the following assertions hold.

\medskip
\begin{compactitem}
\item[\rm(i)] The group $G_{S\cup S_0}^T(k)(p)$ has cohomological dimension~$2$ and the cup-product
\[
H^1(G_{S\cup S_0}^T(k)(p),\F_p) \otimes H^1(G_{S\cup S_0}^T(k)(p),\F_p) \lang H^2(G_{S\cup S_0}^T(k)(p),\F_p)
\]
is surjective. \smallskip
\item[\rm (ii)] We have
\[
k_{S\cup S_0}^T(p)_\p=k_\p(p)
\]
for all\/ $\p\in S\cup S_0$, i.e.\  the extension of global fields $k_{S\cup S_0}^T(p)|k$ realizes the maximal $p$-extension $k_\p(p)$ of the local field~$k_\p$ for all $\p\in S\cup S_0$.

\smallskip
\item[\rm (iii)] The decomposition groups of the primes in $S$ establish a free product in the group $G_{S\cup S_0}^T(k)(p)$, i.e.\ the natural homomorphism
    \[
    \freeproductmed_{\p\in S(k_{S_0}^{S\cup T}(p))} \Gal(k_\p(p)|k_\p) \lang \Gal(k_{S\cup S_0}^T(p)|k_{S_0}^{S\cup T}(p))
    \]
    is an  isomorphism of pro-$p$-groups.

\smallskip
\item[\rm (iv)] For each  discrete $p$-torsion $G_{S\cup S_0}^T(p)$-module $M$, the edge morphisms of the Hochschild-Serre spectral sequence for the universal pro-$p$-covering
    \[
    H^i(G_{S\cup S_0}^T(k)(p), M) \lang H^i_\et (X\sm (S\cup S_0), T, M)
    \]
    are isomorphisms for all $i\geq 0$. Here $X$ is the uniquely defined one-dimensional, regular and connected scheme which is proper over $\Spec(\Z)$ and has function field $k$, and $H^i_\et (X\sm (S\cup S_0), T, M)$ denotes the  \'{e}tale cohomology of the marked arithmetic curve $(X\sm (S\cup S_0),T)$  with values in the sheaf defined by $M$ (see section~\ref{mark-sec}).
\end{compactitem}

\end{theorem}

\begin{remarks} 1. The reason why we excluded the prime number $p=2$ in Theorem~\ref{haupt} is not the usual problem with the real places, but that there is no good theory of mild pro-$2$-groups at hand at the moment. 

\noindent
2. Naturally the question arises whether properties (i)--(iv) then also hold for the group  $G_{S\cup S_0'}^T(k)(p)$, where $S_0'$ is an arbitrary finite set of places which contains $S_0$ and is disjoint from $T$. We will prove this in section~\ref{erwsec} under the assumption that no prime in $S_0'\sm S_0$  splits completely in the extension $k_{S\cup S_0}^T(p)|k$ (see Theorem~\ref{enlarge}).

\noindent
3. Similar results for the full group $G_S^T(k)$, i.e.\ without passage to the maximal pro-$p$ factor group, seem to be out of reach at the moment.
\end{remarks}

In the number field case, the places above $p$ do not play a special role in Theorem~\ref{haupt}, in particular, we do not assume that $S$ contains the set $S_p$ of primes dividing~$p$.  But even in the number field case with $S\supset S_p$ and $T=\varnothing$, Theorem~\ref{haupt} provides new information: assertion (ii) was known up to now only if $k$ contains a primitive $p$-th root of unity (Kuz'min's Theorem, see \cite{kuz} or \cite{NSW},  10.8.4), as well as for certain CM-fields (see \cite{muk} or \cite{NSW}, X \S8 Exercise). By (iii), after adding finitely many primes to $S\supset S_p$, the decomposition groups of the primes above~$p$ establish a free product inside the group $G_S(k)(p)$. So far this was only known for sets of places $S$ of Dirichlet density~$1$ (see \cite{NSW}, 9.4.4), but not for finite sets of places.

In the ``tame'' case $S\cap S_p=\varnothing$ with $T=\varnothing$, assertions (i), (ii) and (iv) were shown in \cite{kpi1}, however only if  $p\nmid \# \Cl(k)$ and $\zeta_p\notin k$. Previous results in the  ``mixed'' case $\varnothing \varsubsetneq S\cap S_p \varsubsetneq S_p$ had been achieved by  K.~Wingberg \cite{wing}, Ch.~Maire \cite{maire} and D.~Vogel \cite{Vo}.

\smallskip
We will see in section~\ref{dualsec} that Theorem~\ref{haupt} provides a large class of examples of  $G_S^T(k)(p)$ being a pro-$p$ duality group (see Theorem~\ref{dualmod}). In the number field case with $S\supset S_p$ and $T=\varnothing$, this was known  if $k$ contains a  primitive $p$-th root of unity by results of Wingberg (see \cite{wi} or \cite{NSW}, 10.9.8). In the case  $\zeta_p\notin k$, there existed merely results for real-abelian fields and for certain CM-fields  (see \cite{NSW}, 10.9.15 and the Remark below).

\smallskip
Essential for the proof of Theorem~\ref{haupt} are arithmetic duality theorems, Hasse principles for the cohomology, Labute's theory of mild pro-$p$-groups and the technique of free products of bundles of profinite groups as developed in Chapter~IV of~\cite{NSW}.

\smallskip
The author thanks Ph.~Lebacque for his comments about an earlier version of this paper.

\section{The \'{e}tale site of a marked curve} \label{mark-sec}

Let  $Y$ be a one-dimensional, noetherian and regular scheme and let $T$ be a finite set of closed points of~$Y$. As usual, we denote by  $\Et(Y)$ the category of \'{e}tale morphisms of finite type  $Y'\to Y$.

\begin{definition}
The category $\Et(Y,T)$ is the full subcategory of $\Et(Y)$ which consists of all objects  $f: Y'\to Y$ such that for each closed point $y'\in Y'$ with $y=f(y')\in T$ the residue field extension $k(y')|k(y)$ is trivial. The \'{e}tale site $(Y,T)_\et$ of the curve $Y$ marked in $T$ consists of the category $\Et(Y,T)$ with surjective families as coverings.
\end{definition}

Obviously, we have $(Y,\varnothing)_\et=Y_\et$. For $T_1\subset T_2$, we have a natural morphism $\iota: (Y,T_1)_\et\to (Y,T_2)_\et$ and hence homomorphisms
\[
H^i_\et(Y,T_2, F) \lang H^i_\et(Y,T_1,\iota^*F)
\]
for each sheaf $F$ of abelian groups on $(Y,T_2)_\et$  and for all $i\geq 0$. For a proper closed (hence finite) subset $M\subset Y$, we define the local cohomology groups $H^*_M(Y,T,-)$ in the usual manner as the right derivatives of the functor
\[
F \longmapsto \ker \big(\Gamma(Y,T,F) \to \Gamma (Y\sm M, T\sm M, F)\big).
\]
In the same way as for usual the \'{e}tale cohomology, one proves excision, i.e.\ we have
\[
H^i_M(Y,T, F)\cong \bigoplus_{x\in M} H^i_x(Y_x^h, T_x^h, F),
\]
where $Y_x^h$ is the  henselization of $Y$ in $x$ and $T_x^h$ is the preimage of $T$ in $Y_x^h$ (i.e.\ $T_x^h=\{x\}$ if $x\in T$ and  $T_x^h=\varnothing$ otherwise). Also the construction of the cup-product is completely analogous to the usual \'{e}tale site:  for sheaves \pagebreak $F_1,F_2$ on $(Y,T)_\et$ and $i,j\geq 0$, we have a cup-product pairing
\[
H^i_\et(Y, T, F_1) \times H^j_\et(Y, T, F_2) \stackrel{\cup}{\lang} H^{i+j}_\et(Y, T, F_1\otimes F_2),
\]
satisfying the using properties.

\bigskip
Also in an analogous way, one constructs the fundamental group. We consider the full subcategory $\FEt(Y, T)$ of the finite morphisms $Y' \to Y$ in $\Et(Y, T)$. This category satisfies the axioms of a Galois category (\cite{sga1}, V,~4). Choosing a geometric point $\bar x$ in $Y \sm T$, we have the fibre functor
\[
\FEt(Y, T) \lang (\textit{Sets}),\  (Y'\to Y)\mapsto \Mor_{Y}(\bar x, Y'),
\]
whose group of automorphisms is the \'{e}tale fundamental group of $(Y,T)$ by definition. We denote this group by  $\pi_1^\et(Y,T,\bar x)$. Assume that $Y$ is connected. Then the fundamental groups to different base points are isomorphic, the isomorphism being canonical up to inner automorphisms, and we will frequently exclude the base point from the notation. The fundamental group is profinite and classifies \'{e}tale coverings of $Y$ in which every point of $T$ splits completely.
The group $H^1_\et(Y, T,\F_p)$ classifies cyclic coverings of degree~$p$ of $(Y, T)$. Therefore we have isomorphisms
\[
H^1 (\pi_1^\et(Y,T)(p), \F_p) \liso H^1 (\pi_1^\et(Y,T), \F_p) \liso H^1_\et(Y, T,\F_p),
\]
where $\pi_1^\et(Y,T)(p)$ denotes the maximal pro-$p$ factor group of~$\pi_1^\et(Y,T)$. This isomorphism is visible in the Hochschild-Serre spectral sequence as follows. We consider the universal pro-$p$-covering  $\widetilde{(Y,T)}(p)$ of $(Y,T)$. This is a pro-object in $\FEt(Y, T)$ and the projection
\[
\widetilde{(Y,T)}(p) \lang (Y, T)
\]
is Galois with group  $\pi_1^\et(Y,T)(p)$. For a discrete $p$-torsion $\pi_1^\et(Y,T)(p)$-module $M$, we obtain the spectral sequence
\[
E_2^{ij}= H^i\Big(\pi_1^\et(Y,T)(p), H^j_\et \big(\widetilde{(Y,T)}(p), M\big) \Big)  \Rightarrow H^{i+j}_\et(Y, T, M)
\]
and, in particular, edge morphisms
\[
\phi_{i,M}: H^i(\pi_1^\et(Y,T)(p), M) \longrightarrow H^i_\et(Y, T, M), \quad i\geq 0.
\]
\begin{lemma}\label{kpi1lem} Let  $Y$ be a connected, one-dimensional, noetherian and regular scheme and let $T$ be a finite set of closed points of\/~$Y$. Then, for every discrete $p$-torsion $\pi_1^\et(Y,T)(p)$-module $M$, the homomorphisms $\phi_{0,M}$ and $\phi_{1,M}$ are isomorphisms and $\phi_{2,M}$ is injective. Furthermore, the following assertions are equivalent.

\medskip
\begin{compactitem}
\item[\rm (i)] $\phi_{i,M}$ is an isomorphism for all $i$ and all $M$. \smallskip
\item[\rm (ii)] $\phi_{i,\F_p}$ is an isomorphism for all $i$. \smallskip
\item[\rm (iii)] $H^i_\et(\widetilde{(Y,T)}(p), \F_p)=0$ for all $i\geq 1$.
\end{compactitem}
\end{lemma}

\begin{proof}
By construction, we have $H^1_\et(\widetilde{(Y,T)}(p), \F_p)=0$, showing the first assertion. As $\pi_1^\et(Y,T)(p)$ is a pro-$p$-group, $\F_p$ is the only simple $p$-torsion $\pi_1^\et(Y,T)(p)$-module. Hence (ii) implies (i)  for finite $M$. Since cohomology of profinite groups as well as \'{e}tale cohomology commute with filtered direct limits, (ii) implies (i) for an arbitrary $p$-torsion $\pi_1^\et(Y,T)(p)$-module $M$.  We can read off the implication  (iii)$\Rightarrow$(ii) directly from the spectral sequence. Finally, assume that (i) holds. Each class in $H^i_\et(\widetilde{(Y,T)}(p), \F_p)$ is already in $H^i_\et(Y',T', \F_p)$ for a finite intermediate covering $(Y',T')$ of $\widetilde{(Y,T)}(p)\to(Y,T)$. The covering $(Y',T')$ corresponds to an open subgroup  $U\subset \pi_1^\et(Y,T)(p)$. Applying  (i) to the $\pi_1^\et(Y,T)(p)$-module $\mathrm{Ind}^U_{\pi_1^\et(Y,T)(p)} \F_p$, we see that each $\alpha \in H^i_\et(Y', T', \F_p)$ vanishes in $H^i_\et(Y'', T'', \F_p)$ for a suitable chosen finite intermediate covering $(Y'',T'')$ of $\widetilde{(Y,T)}(p)\to(Y',T')$. This shows (iii) and finishes the proof.
\end{proof}

\begin{definition} \label{kpi1def} If the equivalent assertions of  Lemma~\ref{kpi1lem} are satisfied, we say that $(Y, T)$ has the {\bf\boldmath $K(\pi,1)$-property for $p$}.
\end{definition}

\begin{remark}
Lemma~\ref{kpi1lem} and Definition~\ref{kpi1def} naturally extend to pro-objects. A marked curve $(Y,T)$ has the $K(\pi,1)$-property for $p$ if and only if this is the case for its universal pro-$p$-covering $\widetilde{(Y,T)}(p)$.
\end{remark}

\section{Calculations of cohomology groups}

In the following, we want to calculate the \'{e}tale cohomology of marked arithmetic curves. Let $k$ be a local or global field and let $p\neq \text{char}(k)$ be a fixed prime number. We denote the group of $p$-th roots of unity by $\mu_p$ and put $\delta=1$ if $\mu_p\subset k$ and $\delta=0$ otherwise. All cohomology groups take values in the constant sheaf $\F_p$, which we exclude from the notation. Furthermore, we put
\[
h^i(-)=\dim_{\F_p} H^i_\et(-)\quad(=\dim_{\F_p} H^i_\et(-,\F_p)\;)\,.
\]
We start with a local computation. For a local field $k$ which is not an extension of $\Q_p$, we use the convention $[k:\Q_p]=0$. Furthermore, we denote the unramified cohomology by  $H^i_\nr(k)$ and set
\[
H^i_{/\nr}(k)= H^i(k)/H^i_\nr(k).
\]
By $A^\vee$ we denote the Pontryagin dual of  $A$.
\renewcommand{\arraystretch}{1.2}
\begin{proposition}\label{localcoh} Let $k$ be a non-archimedean local field with $\text{char}(k)\neq p$.
Let $X=\Spec(\O_k)$, $x$ the closed point of~$X$ and $T=\varnothing$ or~$T=\{x\}$.
Then the local \'{e}tale cohomology groups $H^i_x(X,T)$ vanish for $i\leq 1$ and $i\geq 4$. Moreover, we have
\[
H^2_x(X,T)\cong \left\{
\begin{array} {ll}
H^1_{/\nr}(k)& \text{ if } T=\varnothing,\\
H^1(k) & \text{ if } T=\{x\},
\end{array}\right.
\]
hence
\[
h^2_x(X,T)= \delta+[k:\Q_p] + \# T \,.
\]
Furthermore, we have $H^3_x(X,T)\cong H^2(k)\cong \mu_p(k)^\vee$, hence $h^3_x(X,T)=\delta$, and the following Euler-Poincar\'{e} characteristic formula holds:
\[
\sum_{i=0}^3 (-1)^i h^i_x(X,T)= [k:\Q_p] + \# T\,.
\]
\end{proposition}
\renewcommand{\arraystretch}{1}

\begin{proof}
Since $X$ is henselian, we have isomorphisms  $H^i_\et(X) \cong H^i_\et(x)$ for all~$i$. Furthermore, there are no nontrivial \'{e}tale covers of $(X,T)$ if $T=\{x\}$.   Therefore we obtain
\[
h^i(X,T)= \left\{
\begin{array}{cl}
1 & \hbox{ f\"{u}r } i=0,\\
1-\# T & \hbox{ f\"{u}r } i=1,\\
0 & \hbox{ f\"{u}r } i\geq 2.
\end{array}
\right.
\]
Furthermore, $X\sm \{x\}=\Spec(k)$, hence $H^i_\et(X\sm \{x\})\cong H^i(k)$. The local duality theorem (see \cite{NSW}, Theorem 7.2.6) shows $H_\et^2(X\sm \{ x\})\cong \mu_p(k)^\vee$ and, by \cite{NSW}, Corollary 7.3.9, we obtain
\[
h^1(X\sm \{ x\})= 1+\delta+[k:\Q_p].
\]
Finally, the natural homomorphism $H^1_\et(X) \to H^1_\et(X\sm \{ x\})$  is injective. Therefore the statement of the proposition follows from the exact excision sequence
\[
\cdots \to H^i_x(X,T) \to H^i_\et(X,T) \to H^i_\et(X\sm \{x\}) \to H^{i+1}_x(X,T) \to \cdots\;.
\]
\end{proof}

Now let $k$ be a global field and let $X$ be the uniquely defined one-dimensional, regular and connected scheme which is proper over $\Spec(\Z)$ and has function field~$k$ (i.e.\ $X=\Spec(\O_k)$ if $k$ is a number field and $X$ is a smooth projective curve over a finite field if $k$ is a function field). Let $S$ and $T$ be disjoint finite sets of non-archimedean places of $k$, i.e.\ disjoint finite sets of closed points of~$X$.

\medskip\noindent
We denote the set of archimedean places of $k$ by  $S_\infty$ (hence $S_\infty=\varnothing$ if $k$ is a function field).   In Galois terminology (and without mention the base point), we have
\[
\pi_1^\et(X\sm S, T)= G_{S\cup S_\infty}^T(k):=\Gal (k_{S\cup S_\infty}^T|k),
\]
where $k_{S\cup S_\infty}^T$ is the maximal extension of $k$ which is unramified outside $S\cup S_\infty$ and completely split at~$T$. If $K|k$ is an intermediate extension of $k_{S\cup S_\infty}^T|k$, we denote by
\[
(X\sm S, T)_K
\]
the normalization $(X\sm S)_K$ of the curve $X\sm S$ in $K$, marked in the set  $T(K)$ of prolongations of places of $T$ to~$K$. If $K|k$ is finite, then $(X\sm S, T)_K$ is an object in $\FEt(X\sm S,T)$, and a pro-object otherwise.

\medskip
Now let $p\neq\text{char}(k)$ be a fixed prime number. We denote the maximal elementary abelian $p$-extension of $k$ inside $k_{S\cup S_\infty}^T$ by
\[
k_{S\cup S_\infty}^{T,\el}
\]
and note that
\[
G(k_{S\cup S_\infty}^{T,\el}|k) \cong H^1_\et(X\sm S, T)^\vee.
\]
We denote the set of primes dividing $p$ by $S_p$ (hence $S_p=\varnothing$ in the function field case). Assuming $p\neq 2$ or $k$ totally imaginary if $k$ is a number field, we may ignore the archimedean primes, i.e.\ we have
\[
G_{S}^T(k)(p)=G_{S\cup S_\infty}^T(k)(p)\quad \big(= \pi_1^\et(X\sm S, T)(p)\big).
\]
Part (iv) of our main result Theorem~\ref{haupt} says that, for $p\neq 2$, $p\neq \text{char}(k)$, we can achieve the $K(\pi,1)$-property for $(X\sm S,T)$ by adding finitely many primes to~$S$. This is well known, even without adding primes, in the special case $S\supset S_p$ and $T=\varnothing$:

\begin{proposition} \label{pinS}
Let $k$ be a global field,  $p\neq \text{char}(k)$ a prime number and  $S\supset S_p$ a finite, nonempty set of places of\/ $k$. Then, for every discrete $p$-torsion  $G_{S\cup S_\infty}(k)$-module $M$, the natural maps
\[
H^i(G_{S\cup S_\infty}(k), M) \lang H^i_\et(X\sm S, M)
\]
are isomorphisms for all\/ $i\geq 0$. If, moreover, $M$ is a discrete $p$-torsion $G_{S\cup S_\infty}(k)(p)$-module, then also the natural maps
\[
H^i(G_{S\cup S_\infty}(k)(p), M) \lang H^i(G_{S\cup S_\infty}(k), M)
\]
are isomorphisms for all\/ $i\geq 0$.
\end{proposition}

\begin{proof}
See \cite{Mi}, II, Proposition 2.9, for a proof of the first statement.
The second statement is a result of  O.~Neumann, see \cite{NSW}, Corollary 10.4.8.
\end{proof}

From now on we use the following notation, where $p$ is a prime number different from $\mathrm{char}(k)$.
\begin{tabbing}
\quad \= $r_1$\hspace*{1cm} \= the number of real places of\/ $k$  \\
\>$r_2$\ \> the number of complex places of\/ $k$\\
\>$r$\ \> $=r_1+r_2$, the number of archimedean places of\/ $k$\\
\> $\delta$ \> equals $1$ if $\mu_p\subset k$ and $0$ otherwise\\
\> $\delta_\p$ \> equals $1$ if $\mu_p\subset k_\p$ and $0$ otherwise\\
\> $\Cl_S(k)$\> $=\Pic(X\sm S)$, the $S$-ideal class group of\/ $k$\\
\> $E_{k,S}$ \> $=H^0(X\sm S,\G_m)$, the group of $S$-units of\/ $k$\\
\> $_n A$\> $=\ker(A \stackrel{\cdot n}{\to} A)$, where $A$ is an abelian group and $n\in \N$ \\
\> $A/n$\> $=\coker(A \stackrel{\cdot n}{\to} A)$, where $A$ is an abelian group and  $n\in \N$.\\
\end{tabbing}
If $k$ has positive characteristic, then $r_1=r_2=r=0$, and $E_{k,\varnothing}$ is the multiplicative group of the finite field of constants. As before, we exclude the constant coefficients $\F_p$ from the notation in cohomology groups. We make the assumption $p\neq 2$ or $k$ totally imaginary in the number field case everywhere in the paper. Therefore we make the following notational convention:

\medskip
\begin{minipage}{11cm}
{\em `Set of places' always means `set of non-archimedean places'. }
\end{minipage}

\medskip\noindent

\begin{proposition}\label{globalchi} Let $S$ and\/ $T$ be disjoint, finite sets of places of the global field~$k$. We assume $p\neq 2$ or $k$ totally imaginary if\/ $k$ is a number field. Then the groups $H^i_\et(X\sm S,T)=0$ vanish for\/ $i\geq 4$ and we have
\[
\chi(X\sm S,T):= \sum_{i=0}^3 (-1)^i h^i(X\sm S,T)= r + \# T - \sum_{\p \in S\cap S_p} [k_\p:\Q_p]\;.
\]
\end{proposition}

\begin{proof} The assertion is well known in the case $T=\varnothing$ and $S\supset S_p$: by \cite{Mi}, II, Theorem~2.13\,(a), we have $\chi(X\sm S)=-r_2$, and furthermore
\[
r - \sum_{\p\in S_p} [k_\p :\Q_p] = r_1+r_2 - [k:\Q]= -r_2,
\]
where we set $[k:\Q]=0$, if $k$ is a function field. We consider the exact excision sequence
\[
\cdots \to \bigoplus_{\p \in S\cup T} H^i_\p(X_\p, T_\p) \to H^i_\et(X,T) \to H^i_\et (X\sm (S\cup T)) \to  \cdots \ ,
\]
where $X_\p$ denotes the completion of $X$ at $\p$ (the cohomology groups of the completion and the henselization coincide). We obtain the result for  $S=\varnothing$ and arbitrary $T$ by using this excision sequence for $S=S_p$ and by applying Proposition~\ref{localcoh}.  Finally, we obtain the result for arbitrary $S$ from the case $S=\varnothing$ by using the excision sequence
\[
\cdots \to \bigoplus_{\p \in S} H^i_\p(X_\p) \to H^i_\et(X,T) \to H^i_\et (X\sm S,  T) \to  \cdots \ ,
\]
and another application of Proposition~\ref{localcoh}.
\end{proof}

\begin{corollary}\label{ausschneid2} Let $S$ and\/ $T$ be disjoint, finite sets of places of the global field~$k$. Then we have an exact sequence
\[
 H^1_\et(X\sm S,T) \hookrightarrow H^1_\et(X\sm S) \to \bigoplus_{\p\in T} H^1_{\nr}(k_\p) \to H^2_\et(X\sm S,T) \twoheadrightarrow H^2_\et(X\sm S)
\]
and isomorphisms $H^i(X\sm S,T)\liso H^i(X\sm S)$ for $i\geq 3$.
\end{corollary}

\begin{proof}
This follows from Proposition~\ref{localcoh} by comparing the excision sequences for  $X\sm S$ and $X\sm (S\cup T)$ as well as for $(X\sm S,T)$ and $X\sm (S\cup T)$.
\end{proof}

To obtain formulae for the individual cohomology groups, we introduce the Kummer group $V_S^T(k)$.
Let $S$ and $T$ be disjoint, finite sets of places of $k$ and let $p\neq \textrm{char}(k)$ be a fixed prime number, which we exclude from the notation. We put
\smallskip\noindent
\[
V^T_{S}(k):=
\{ a \in k^\times\,|\, a \in k_v^{\times p}
\hbox{ for } v \in {S}\,\,\hbox{and}\,\, a \in U_v k_v^{\times p}
\hbox{ for } v \notin T \}/ k^{\times p}, \smallskip
\]
where $U_v$ denote the unit group of the local field  $k_v$  (convention:
$U_v=k_v^\times$, if\/~$v$ is archimedean). In terms of flat cohomology, we have
\[
V_S^T(k)=\ker\big(H^1_\fl(X\sm (S\cup T),\mu_p) \lang \prod_{\p\in S} H^1(k_\p,\mu_p)\big).
\]

\begin{lemma}\label{VSchange}
For $S=\varnothing$, there is a natural exact sequence
\[
0\longrightarrow E_{k,T} /p \longrightarrow V_\varnothing^T(k) \longrightarrow \null_p \Cl_T(k) \longrightarrow 0\,.
\]
In particular, we have
$
\dim_{\F_p} V_\varnothing^T(k) =  \dim_{\F_p} \null_p \Cl_T(k) + r-1+\# T + \delta
$, unless $k$ is a function field and $T=\varnothing$, where $\dim_{\F_p} V_\varnothing^\varnothing(k)=\dim_{\F_p} \null_p \Cl(k)+\delta$.

\smallskip
For arbitrary $S$ and any additional place  $\p\notin S\cup T$, we have an exact sequence
\[
0 \longrightarrow V_{S \cup \{\p\}}^T(k)\stackrel{\phi}{\longrightarrow} V_S^T(k)  \longrightarrow U_\p k_\p^{\times p}/k_\p^{\times p}.
\]
In particular, $V_S^T(k)$ is finite. For $\p\notin S_p$, we have
$
\dim_{\F_p} \coker(\phi)\leq \delta_\p$.
\end{lemma}

\begin{proof}
Sending $a\in V_\varnothing^T(k)$ to the class $[\a]$ of the uniquely defined fractional $T$-Ideal $\a$ with $(a)={\mathfrak a}^p$, defines a surjective homomorphism $V_\varnothing^T(k) \to \null_p\Cl_T(k)$ whose kernel is $E_{k,T}/p$. This, together with Dirichlet's unity theorem shows the first statement. The second exact sequence follows directly from the definitions of the objects occurring. For $\p\notin S_p$, we have $\dim_{\F_p} U_\p k_\p^{\times p}/k_\p^{\times p} = \delta_\p$, showing the last statement.
\end{proof}

\begin{theorem} \label{globcoh} Let $S$ and\/ $T$ be disjoint, finite sets of places of the global field~$k$. Assume $p\neq 2$ or $k$ totally imaginary if $k$ is a number field.  Then $H^i_\et(X\sm S,T)=0$ for $i\geq 4$
and
\[
\renewcommand{\arraystretch}{1.3}
\begin{array}{lcl}
h^0(X \sm S,T)&=& 1\, ,\\
h^1(X \sm S,T) & = & 1+  \ds\sum_{\p\in S} \delta_\p
               - \delta + \dim_{\F_p}V_S^T(k) +  \ds\sum_{\p\in S\cap S_p} [k_\p:\Q_p] -r - \# T,\\
h^2(X \sm S,T)& =& \ds\sum_{\p\in S} \delta_\p
               - \delta + \dim_{\F_p}V_S^T(k)  +\theta \, ,\\
h^3(X \sm S,T)&=& \theta\, .\\
\end{array}
\renewcommand{\arraystretch}{1}
\]
Here $\theta=1$  if\/ $\delta =1$ and $S=\varnothing$, and $\theta=0$ otherwise.
\end{theorem}

\begin{proof}
The statement about $h^0$ is trivial and the vanishing of the cohomology in dimension greater than or equal to~$4$ is already contained in Proposition~\ref{globalchi}.  Artin-Verdier duality (see \cite{Ma}, 2.4 or \cite{Mi}, Theorem~3.1) or \'{e}tale Poincar\'{e} duality (\cite{Mi1}, V, Corollary~2.3), respectively, show that
\[
H^3_\et(X)^\vee\cong\Hom_X(\F_p, \G_m)=\mu_p(k).
\]
Furthermore, the vanishing of $H_\et^i(X_\p,T_\p)$ for $i\geq 2$ and the local duality theorem (\cite{NSW}, Theorem 7.2.6) show isomorphisms $H^3_\p(X_\p,T_\p)^\vee\cong H^2(k_\p, \F_p)^\vee\cong \mu_p(k_\p)$ for every point  $\p$ of~$X$. By Corollary~\ref{ausschneid2}, we have an isomorphism $H^3_\et(X,T) \stackrel{\sim}{\to} H^3_\et(X)\cong \mu_p(k)^\vee$. We consider the excision sequence
\[
\bigoplus_{\p\in S} H^3_\p(X) \stackrel{\alpha}{\to} H^3_\et (X,T) \stackrel{\beta}{\to} H^3_\et(X\sm S,T) \to \bigoplus_{\p\in S} H^4_\p(X_\p).
\]
The right hand term vanishes by Proposition~\ref{localcoh}, hence $\beta$ is surjective.  The dual map to $\alpha$ is the map
\[
 \mu_p(k) \to \bigoplus_{\p\in S} \mu_p(k_\p),
\]
which is injective, unless  $\delta=1$ and $S=\varnothing$.
We obtain $h^3(X \sm S,T)=1$ if $\delta=1$ and $S = \varnothing$, and $0$ otherwise.
Using the isomorphism $H^1(G_S^T(k)) \stackrel{\sim}{\rightarrow} H^1_\et(X\sm S,T)$, we obtain the assertion about  $h^1$  from the calculation of the first cohomology of the group $G_S^T(k)$, see  \cite{NSW}, Theorem 10.7.10. Finally, the result for $h^2$ follows by using the Euler-Poincar\'{e} characteristic formula in Proposition~\ref{globalchi}.
\end{proof}

The vanishing of the cohomology in dimension greater than $2$ and Lemma~\ref{kpi1lem} imply the following
\begin{corollary}
\label{kpi1crit} For $S\neq \varnothing$ and each intermediate extension $K|k$ of $k_S^T(p)|k$, the following assertions are equivalent:

\smallskip
\begin{compactitem}
\item[\rm (i)] $(X\sm S,T)_K$ has the $K(\pi,1)$-property for $p$.\smallskip
\item[\rm (ii)] The homomorphism $\phi_{2,\F_p}: H^2(G_S^T(K)(p)) \to H^2_\et((X\sm S,T)_K) $ is surjective and  $\cd\, G_S^T(K)(p)\leq 2$.
\end{compactitem}
\end{corollary}

\section{Duality for marked arithmetic curves}

Next we prove a duality theorem for marked arithmetic curves. For a sheaf $F$ on $(X\sm S, T)_\et$, we consider the  Shafarevich-Tate groups
\[
\Sha^i(k,S,T, F)= \ker\big(H^i_\et(X\sm S, T, F) \lang \bigoplus_{\p\in S} H^i(k_\p, F) \big).
\]
We omit the constant coefficients $F=\F_p$ from the notation and we also omit~$T$ if $T=\varnothing$. Furthermore, we put
\[
\Be_S^T(k)= V_S^T(k)^\vee.
\]

\begin{theorem}  \label{bstdual} Let  $k$ be a global field and let $p$ be a prime number different from $\mathrm{char}(k)$. Assume $p\neq 2$ or $k$ totally imaginary if $k$ is a number field. Let  $S$ and $T$ be disjoint, finite sets of places of~$k$. Then there is a natural isomorphism
\[
\Sha^2(k,S,T) \liso \Be_S^T(k).
\]
In particular, $\Sha^2(k,S,T)$ is finite.
\end{theorem}

\begin{remark} By passing to the limit over all finite subsets, Theorem~\ref{bstdual} generalizes to the case of an arbitrary, i.e.\ not necessarily finite set of places~$S$.
\end{remark}

\begin{proof}[Proof of Theorem~\ref{bstdual}] We choose a sufficiently large finite set of places $\Sigma$ such that $S\cup T \cup S_p\subset \Sigma$ and   $\Pic\big((X\sm \Sigma)_{k(\mu_p)}\big)$ is $p$-torsion free. Then the groups $\Sha^1(k,\Sigma)$ and $\Sha^1(k, \Sigma, \mu_p)$ vanish.
From this follows the vanishing of $\Sha^2(k,\Sigma)\cong \Sha^1(k,\Sigma,\mu_p)^\vee$ (by Proposition~\ref{pinS} and Poitou-Tate duality, \cite{NSW}, Theorem 8.6.7) and the exactness of the sequence
\[
0 \lang V_S^T(k) \lang H^1_\et(X\sm \Sigma,\mu_p) \lang \prod_{\p\in S} H^1(k_\p,\mu_p) \times \back \prod_{\p\in \Sigma\backslash (S \cup T)} \back \!\!\!  H^1_{/\nr}(k_\p,\mu_p)\;.
\]
We dualize this sequence  and consider the commutative diagram
\[
\xymatrix@C=.5cm{&\ds\prod_{\p\in S} H^1(k_\p) \times \back \prod_{\p\in \Sigma\backslash (S \cup T)} \back   \!\!\! H^1_{\nr}(k_\p)\ar[r]\ar@{^{(}->}[d]&H^1_\et(X\sm \Sigma,\mu_p)^\vee\ar@{=}[d]\ar@{->>}[r]&\Be_S^T(k)\\
H^1_\et(X\sm \Sigma)\ar@{^{(}->}[r] & \ds\prod_{\p\in \Sigma} H^1(k_\p)\ar@{->>}[r]\ar@{->>}[d]&H^1_\et(X\sm \Sigma,\mu_p)^\vee\\
&\ds\prod_{\p\in T} H^1(k_\p) \times \back \prod_{\p\in \Sigma\backslash (S \cup T)} \back \!\!\!  H^1_{/\nr}(k_\p)\;.
}
\]
The middle row is part of the long exact sequence of Poitou-Tate (see \cite{NSW}, Theorem 8.6.10), hence exact. The snake lemma yields the exactness of
\[
H^1_\et(X\sm \Sigma) \lang \prod_{\p\in S} H^1(k_\p) \times \back \prod_{\p\in \Sigma\backslash (S \cup T)} \back   \!\!\! H^1_{/\nr}(k_\p) \lang \Be_S^T(k) \lang 0\,. \leqno \mathrm{(I)}
\]
Using the commutative diagram
\[
\xymatrix{H^2_\et(X\sm S, T)\ar[d]\ar[r]&H^2_\et(X\sm \Sigma)\ar@{^{(}->}[d]\\
\ds\prod_{\p\in S} H^2(k_\p)\ar@{^{(}->}[r] & \ds\prod_{\p\in \Sigma} H^2(k_\p),}
\]
we obtain
\[
\Sha^2(k,S,T)=\ker\big( H^2_\et(X\sm S, T)\lang H^2_\et(X\sm \Sigma)\big).
\]
Therefore excision yields the exact sequence
\[
H^1_\et(X\sm \Sigma) \lang \prod_{\p\in S} H^1(k_\p) \times \back \prod_{\p\in \Sigma\backslash (S \cup T)} \back   \!\!\! H^1_{/\nr}(k_\p) \lang \Sha^2(k,S,T) \lang 0\, , \leqno \mathrm{(II)}
\]
and a comparison of  (I) and (II) shows the statement of the theorem.
\end{proof}

\begin{corollary} \label{h2lokali} Let $k$ be a global field and let $p$ be a prime number different from  $\mathrm{char}(k)$. Assume $p\neq 2$ or $k$ totally imaginary if $k$ is a number field.
Let $S$ and $T$ be disjoint, finite sets of places of\/ $k$ with  $V_S^T(k)=0$.
Then the natural map
\[
H_\et^2(X\sm S, T) \lang \prod_{\p\in S} H^2(k_\p)
\]
is injective and an isomorphism if $\delta=0$. If\/  $\delta=1$, then, for every  $\p_0 \in S$, the map
\[
H_\et^2(X\sm S, T) \lang \prod_{\p\in S \backslash \{\p_0\}} H^2(k_\p)
\]
is an isomorphism.
\end{corollary}

\begin{proof} For $\p \in S$, we have
\[
H^3_\p(X,T) \cong H^2(k_\p) \cong \mu_p(k_\p)^\vee.
\]
Furthermore, we have $\Sha^2(k,S,T)=0$ by Theorem~\ref{bstdual}. Therefore excision yields the exact sequence
\[
0 \to H_\et^2(X\sm S, T) \to \prod_{\p\in S} H^2(k_\p) \stackrel{\alpha}\to H_\et^3(X,T).
\]
The map dual to $\alpha$ is $\mu_p(k) \to \prod_{\p\in S} \mu_p(k_\p)$, showing the statement of the corollary.
\end{proof}

Finally, we can annihilate $V_S^T(k)$ in the following sense.

\begin{proposition}\label{vskill} Let $T$ and $\M$ be disjoint sets of places, where  $T$ is finite and $\M$ has Dirichlet density $\delta(\M)=0$. Then there exists a finite set $S$ of places $\p$ with  $N(\p)\equiv 1 \bmod p$ such that $S\cap (T\cup \M)=\varnothing$ and
\[
V_{S}^T(k)=0.
\]
\end{proposition}

\begin{proof}  Let $\varOmega$ be the set of all places $\p$ of $k$ with $N(\p)\equiv 1 \bmod p$. Then the set  $\varOmega(k(\mu_p))$ of prolongations of places in $\varOmega$ to $k(\mu_p)$ has Dirichlet density~$1$.  As  $k^\times/k^{\times p}\cong H^1(k,\mu_p)$ and by the  Hasse principle \cite{NSW}, Theorem 9.1.9\,(ii), the homomorphism
\[
k^\times/k^{\times p} \lang \prod_{\p\in \varOmega\backslash (T\cup\M)} k_\p^\times/k_\p^{\times p}
\]
is injective. Since the $\F_p$-subspace $V_\varnothing^T(k)\subset k^\times/k^{\times p}$ is finite dimensional, we find a finite subset  $S\subset \varOmega\sm (T\cup\M)$ with
\[
V_{S}^T(k)=\ker\big( V_\varnothing^T(k) \lang \prod_{\p\in S} k_\p^\times/k_\p^{\times p} \big)=0.
\]
\end{proof}

\section{Local components}

As before, let $k$ be a global field and let $p$ be a prime number different from  $\text{char}(k)$. We assume $p\neq 2$ or $k$ totally imaginary if $k$ is a number field.  Let $T$ be a finite set of places, $T\neq \varnothing$ if $k$ is a function field, with  $\null_p\Cl_T(k)=0$. As $\Cl_T(k)$ is finite, we have  $\Cl_T(k)(p)=0$ and the exact sequence
\[
0\lang E_{k,T} \lang k^\times \mapr{(v_\q)_\q}\bigoplus_{\q\notin T} \Z \lang \Cl_{T}(k) \lang 0
\]
implies the exactness of
\[
0\lang E_{k,T}/p \lang k^\times/k^{\times p}\lang \bigoplus_{\q\notin T} \Z/p\Z \lang 0.
\]

\begin{definition}
For a prime $\p\notin T$, we denote by $s_\p\in k^\times/k^{\times p}$ an element with  $v_\p(s_\p)\equiv 1\bmod p$ and $v_\q(s_\p) \equiv 0 \bmod p$ for all $\q \notin T \cup \{\p\}$. The element $s_\p$ is well-defined up to multiplication by elements of  $E_{k,T}/p$.
\end{definition}

By $k(\sqrt[p]{E_{k,T}})$ we denote the finite Galois extension of  $k$, obtained by adjoining the $p$-th roots of all $T$-units of $k$. If $\delta=0$ (i.e.\ $\mu_p\not \subset k$), this means that also the $p$-th roots of unity are adjoined. For a place  $\p \notin T$, the extension
$k(\sqrt[p]{E_{k,T}},\sqrt[p]{s_\p})$ is independent from the choice of $s_\p$.

\begin{lemma}  \label{dual}
Let $T$ be finite set of places, $T\neq \varnothing$ if $k$ is a function field, with  $\null_p \Cl_T(k)=0$. Furthermore, let $\q\notin T\cup S_p$ be a prime which splits completely in the extension  $k(\sqrt[p]{E_{k,T}})|k$. Then  $k_{\{\q\}}^{T,\el}|k$ is cyclic of order $p$ and $\q$ ramifies in this extension. Moreover, if\/ $\delta=1$, then
\[
k(\sqrt[p]{E_{k,T}}) k_{\{\q\}}^{T,\el} = k(\sqrt[p]{E_{k,T}},\sqrt[p]{s_\q}).
\]
A place $\p\notin T\cup \{\q\}$ of norm $N(\p)\equiv 1 \bmod p$ splits in $k_{\{\q\}}^{T,\el}|k$ if and only if  $\q$ splits in the extension $k(\sqrt[p]{E_{k,T}},\sqrt[p]{s_\p})|k(\sqrt[p]{E_{k,T}})$.
\end{lemma}

\begin{proof}
By assumption, we have $N(\q)\equiv 1 \bmod p$. Since $\null_p \Cl_T(k)=0$, Lemma~\ref{VSchange} yields an isomorphism $E_{k,T}/p\liso V_\varnothing^T(k)$. As $\q$ splits completely in  $k(\sqrt[p]{E_{k,T}})|k$,  the homomorphism $E_{k,T}/p \to k_{\q}^{\times}/k_{\q}^{\times p}$ is the zero map. Therefore we have an isomorphism $V_{\{\q\}}^T(k) \liso V_\varnothing^T(k)$. In particular,
\[
\dim_{\F_p} V_{\{\q\}}^T(k)= \dim_{\F_p} E_{k,T}/p= \# T + r - 1 + \delta,
\]
which implies $h^1(X\sm \{\q\}, T)=1$ by Theorem~\ref{globcoh}. Because of $\Cl_T(k)(p)=0$, the prime $\q$ ramifies in $k_{\{\q\}}^{T,\el}$. Now assume $\delta=1$ and let $k_{\{\q\}}^{T,\el}=k(\sqrt[p]{\alpha})$, $\alpha\in k^\times/k^{\times p}$. After replacing $\alpha$ by  $\alpha^e$ for a suitable chosen $e \in \F_p^\times$, we have
\[
v_{\q}(\alpha) \equiv 1 \bmod p,\ v_\p(\alpha)\equiv 0 \bmod p \text{ for } \p\neq \q.
\]
Hence $\alpha/s_{\q}$ in $E_{k,T}/p$, implying $k(\sqrt[p]{E_{k,T}}) k_{\{\q\}}^{T,\el} \allowbreak = k(\sqrt[p]{E_{k,T}},\sqrt[p]{s_\q})$ by Kummer theory.

Now let  $\p\notin T\cup \{\q\}$ be a prime of norm $N(\p)\equiv 1 \bmod p$. By class field theory,  $\p$ splits in the extension $k_{\{\q\}}^{T,\el}|k$ if and only if there exists an $s\in k^\times/k^{\times p}$ with

\smallskip
\begin{compactitem}
\item $v_{\mathfrak l}(s) \equiv 0 \bmod p$ for ${\mathfrak l} \notin T \cup \{\p,\q\}$,
\item $v_\p(s) \equiv 1 \bmod p$,
\item $s\in k_\q^{\times p}$.
\end{compactitem}

\smallskip\noindent
If such an $s$ exists, then  $s/s_\p \in E_{k,T}/p$, hence $s_\p \in k_\q^{\times p}$, i.e.\ $\q$ splits in  the extension $k(\sqrt[p]{E_{k,T}},\sqrt[p]{s_\p})|k(\sqrt[p]{E_{k,T}})$. Conversely, if $\q$ splits in this extension, then $s=s_\p$ has the required property and  $\p$ splits in $k_{\{\q\}}^{T,\el}|k$.
\end{proof}

For $x\in H_\et^2(X\sm S, T)$ and $\p\in S$, we denote by $x_\p$ the image of  $x$ under the natural homomorphism  $H_\et^2(X\sm S, T) \lang H^2(k_\p)$ and call it the $\p$-component of~$x$.  For a subset $S'\subset S$, we use the natural embedding to consider elements of $H_\et^1(X\sm S',T)$ also as elements of $H_\et^1(X\sm S, T)$. By $\T_\p \subset G(k_S^{T,\el}|k)$, we denote the inertia group of $\p\in S$. We call  $\chi\in H_\et^1 (X\sm S, T)$ unramified at $\p$ if $\chi(\T_\p)=0$, and ramified otherwise.

\pagebreak
\begin{proposition} \label{komponenten} Let $T$ be finite set of places, $T\neq \varnothing$ if $k$ is a function field, with  $\null_p \Cl_T(k)=0$.
Let $\q\notin T\cup S_p$ be a prime  which splits completely in the extension $k(\sqrt[p]{E_{k,T}})|k$ and let  $\chi_\q$ be a generator of the cyclic (by Lemma~\ref{dual}) group $H_\et^1(X\sm \{ \q \}, T)$.
Furthermore, let $S$ be finite set of primes $\p$ of norm $N(\p)\equiv 1 \bmod p$ with  $S\cap (T\cup \{\q\})=\varnothing$. Then, for an arbitrary element  $\chi \in H_\et^1(X\sm S,T)$, the following holds for the local components of\/ $\chi\cup \chi_\q \in H_\et^2(X\sm (S\cup \{\q\}), T)$:
\[
(\chi \cup \chi_\q)_\q = \left\{
\begin{array}{cl}
0 & \text{if } \chi(\Frob_\q)=0,\\
\neq 0 & \text{ otherwise.}
\end{array}
\right.
\]
For $\p\in S$, we have
\[
(\chi \cup \chi_\q)_\p = \left\{
\begin{array}{cl}
0 & \text{if } \chi \text{ is unramified at $\p$ or $\q$ splits completely} \\
& \text{in } k(\sqrt[p]{E_{k,T}},\sqrt[p]{s_\p})|k(\sqrt[p]{E_{k,T}}), \\
\neq 0 & \text{otherwise.}
\end{array}
\right.
\]
\end{proposition}

\begin{proof}
For a prime $\p$ with  $N(\p)\equiv 1 \bmod p$, the $\F_p$-vector space $H^1(k_\p)$ is two-dimensional, $H^2(k_\p)$ is one-dimensional and the pairing $H^1(k_\p) \times H^1(k_\p) \stackrel{\cup}{\to} H^2(k_\p)$ is perfect. Furthermore,  $\chi_1\cup \chi_2=0$ for $\chi_1,\chi_2 \in H^1_\nr(k_\p)$.

Since $\chi_\q$ ramifies at $\q$ and $\chi$ is unramified at $\q$,  $(\chi \cup \chi_\q)_\q$ is nonzero if and only if the image of  $\chi$ in $H^1_\nr(k_\q)$ is nontrivial, i.e.\ if $\chi(\Frob_\q)\neq 0$.

Now let $\p\in S$. Since $\chi_\q$ is unramified at  $\p$,  $(\chi\cup\chi_\q)_\p$ is nonzero if and only if the image of $\chi_\q$ in $H^1_\nr(k_\p)$ is nonzero and $\chi$ ramifies at $\p$. The first condition is equivalent to $\chi_\q(\Frob_\p)\neq 0$ and, by Lemma~\ref{dual}, this holds if and only if $\q$ does not split in the extension $k(\sqrt[p]{E_{k,T}},\sqrt[p]{s_\p})|k(\sqrt[p]{E_{k,T}})$.
\end{proof}

\section{An auxiliary extension}

\begin{theorem} \label{hilf} Let $T$ and  $\M$ be disjoint sets of places of the global field $k$, where $T$ is finite and  $\M$ has Dirichlet density  $\delta(\M)=0$.
Assume $p\neq 2$ and $p\neq \hbox{char}(k)$. Then there exist a finite set of places $T_0$ and a finite, nonempty set of places $S$ consisting  of primes\/ $\p$ with $N(\p)\equiv 1\bmod p$, such that the following assertions hold.

\smallskip
\begin{compactitem}
\item[\rm (i)]  $S \cap (T\cup T_0 \cup \M)=\varnothing$.\smallskip
\item[\rm (ii)] $(X\sm S, T \cup T_0)$ has the $K(\pi,1)$-property for $p$.\smallskip
\item[\rm (iii)] Each $\p\in S$ ramifies in $k_{S}^{T \cup T_0}(p)$.\smallskip
\item[\rm (iv)] $V_S^{T\cup T_0}(k)=0$. \smallskip
\item[\rm (v)] the cup-product $H^1(G_S^{T\cup T_0}(k)(p)) \otimes H^1(G_S^{T\cup T_0}(k)(p)) \to H^2(G_S^{T\cup T_0}(k)(p))$ is surjective.
\end{compactitem}
\end{theorem}

In the proof of Theorem~\ref{hilf}, we will use the following sufficient criterion for cohomological dimension~$2$.

\begin{theorem}[\cite{kpi1}, Theorem 5.5] \label{mildkrit}
Let $p$ be an odd prime number and let $G$ be a finitely presentable pro-$p$-group.
Suppose $H^2(G)\neq 0$ and that there is a direct sum decomposition $H^1(G)\cong U \oplus V$ such that
\begin{itemize}
\item[\rm(i)] the cup-product $V\otimes V \stackrel{\cup}{\lang} H^2(G)$ is trivial, i.e.\ $v_1\cup v_2=0$ for all $v_1,v_2\in V$,
\item[\rm (ii)] the cup-product $U\otimes V \stackrel{\cup}{\lang} H^2(G)$ is surjective.
\end{itemize}
Then $\cd\, G=2$.
\end{theorem}

\begin{remark} 1. For pro-$p$-groups with one defining relation, this result was known for a long time, see \cite{La2}.

\noindent
2. The criterion in Theorem~\ref{mildkrit} yields more than just cohomological dimension~$2$: the given condition is sufficient for the  {\em mildness} of~$G$. This has been deduced in \cite{kpi1}, \S 5, from the results of Labute's paper \cite{La}. By \cite{La}, Theorem 1.2(c),  mild pro-$p$-groups have  cohomological dimension~$2$.
\end{remark}

The rest of this section is devoted to the proof of Theorem~\ref{hilf}. We start with the (easier) case $\delta=0$, i.e.\ $k$ does not contain a primitive $p$-th root of unity.

\begin{lemma}\label{ewnichtda}
Let $\delta=0$ and let $S=\{\p_1,\ldots,\p_n\}$ be a finite set of places  $\p$ with $N(\p)\equiv 1 \bmod p$. We put $s_i=s_{\p_i}$. Then the extensions \[
k(\mu_p,\sqrt[p]{s_1},\ldots, \sqrt[p]{s_n}), \ k(\sqrt[p]{E_{k,T}}) \text{ and } k_S^{T,\el}(\mu_p)
\]
are linearly disjoint over $k(\mu_p)$.
\end{lemma}

\begin{proof}
By construction, the $\F_p$-vector space spanned by $s_1,\ldots, s_n$ in $k^\times/k^{\times p} E_{k,T}$ is $n$-dimensional. Hence, by Kummer theory, the extensions $k(\mu_p,\sqrt[p]{s_1},\ldots, \sqrt[p]{s_n})|k(\mu_p)$ and $k(\sqrt[p]{E_{k,T}})|k(\mu_p)$ are linearly disjoint. Since both lie in the $1$-eigenspace with respect to the cyclotomic character $\chi_{\text{cycl}}: G(k(\mu_p)|k)\to \F_p^\times$, but $k_S^{T,\el}(\mu_p)|k(\mu_p)$ lies in the $0$-eigenspace, the result follows.
\end{proof}

\begin{proof}[Proof of Theorem~\ref{hilf} in the case $\delta=0$]
We choose $T_0$ such that $\null_p \Cl_{T\cup T_0}(k)=0$ and  $T\cup T_0\neq \varnothing$ if  $k$ is a function field. We simplify the notation and replace $T$ by $T\cup T_0$ in what follows. Now we choose a finite set of $S_0$ places $\p$ with $N(\p)\equiv 1 \bmod p$ such that $S_0\cap (T\cup \M)=\varnothing$ and
\[
V_{S_0\backslash \{\p\}}^T(k)=0 \text{ for each } \p \in S_0.
\]
This can be achieved by applying Proposition~\ref{vskill} twice.  We enumerate the elements in $S_0$, i.e.\  $S_0=\{\p_1,\ldots,\p_m\}$, and put $s_i=s_{\p_i}$.  By Theorem~\ref{globcoh}, the inertia groups $\T_i$ of the places  $\p_i$, $i=1,\ldots, m$, are nontrivial and of order~$p$. We enlarge $S_0$ by $m$ further primes in the following way:

We choose prolongations  $\P_1,\ldots,\P_m$ of $\p_1,\ldots,\p_m$ to $k(\mu_p)$ and consider for a prime  $\mQ$ in $k(\mu_p)$ and  $a\in \{1,\ldots, m\}$ the following condition  $(B_a)$:

\medskip\noindent
\begin{compactitem}
\item $\mQ \notin T(k(\mu_p)) \cup \M(k(\mu_p))$.\smallskip
\item $\Frob_{\mQ} \notin \T_{\P_a} \subset G(k_{S_0}^{T,\el}(\mu_p)|k(\mu_p))$.\smallskip
\item For all $b\neq a$, $\mQ$ splits in $k(\mu_p,\sqrt[p]{s_b})|k(\mu_p)$.\smallskip
\item $\mQ$ is inert in $k(\mu_p,\sqrt[p]{s_a})|k(\mu_p)$.\smallskip
\item $\mQ$ splits completely in $k(\sqrt[p]{E_{k,T}})|k(\mu_p)$.
\end{compactitem}

\medskip\noindent
As the complete splitting of  $\mQ$ in $k(\sqrt[p]{E_{k,T}})$ is one of the conditions,  $(B_a)$ is independent of the choice of the  $s_i$. By Lemma~\ref{ewnichtda} and the Chebotarev density theorem, we find a prime $\mQ_1$ in $k(\mu_p)$ which satisfies condition $(B_1)$.   Then we put $\q_1=\mQ_1 \cap k$.  By Lemma~\ref{dual}, $k_{\{\q_1\}}^{T,\el}$ is cyclic of order~$p$ and $\q_1$ ramifies in this extension. Now we recursively choose primes $\mQ_2,\ldots,\mQ_m$ in $k(\mu_p)$  and put $\q_a=\mQ_a\cap k$ such that

\medskip
\begin{compactitem}
\item $\mQ_a$ satisfies condition $(B_a)$, and \smallskip
\item For $b<a$,   $\mQ_a$ splits in $k_{\{\q_b\}}^{T,\el}(\mu_p)|k(\mu_p)$ and in $k(\mu_p,\sqrt[p]{s_{\q_b}})|k(\mu_p)$.
\end{compactitem}

\medskip\noindent
This is possible since, according to the choice of $\mQ_1,\ldots,\mQ_{a-1}$ and by Lemma~\ref{ewnichtda}, the extensions
\[
k(\mu_p,\sqrt[p]{s_1},\ldots, \sqrt[p]{s_m}, \sqrt[p]{s_{\q_1}},\ldots, \sqrt[p]{s_{\q_{a-1}}} ), \ k(\sqrt[p]{E_{k,T}}) \text{ and } k_{S_0\cup \{ \q_1,\ldots, \q_{a-1}\}}^{T,\el}(\mu_p)
\]
are linearly disjoint over $k(\mu_p)$.  We put
\[
S=\{\p_1,\ldots,\p_m,\q_1,\ldots,\q_m\}.
\]
We have $h^2(X\sm S,T)=2m$ and, by Corollary~\ref{h2lokali}, the natural map
\[
H_\et^2(X\sm S, T) \lang \prod_{i=1}^m H^2(k_{\p_i}) \oplus \prod_{i=1}^m H^2(k_{\q_i})
\]
is an isomorphism.  Let $\eta_a$ be a generator of $H_\et^1(X\sm \{ \q_a\}, T)$. We consider the $m$-dimensional vector space $V$ spanned by  $\eta_1,\ldots,\eta_m$ in $H_\et^1(X\sm S, T)$. We have
\[
H_\et^1(X\sm S, T) \cong H_\et^1(X\sm S_0, T) \oplus V.
\]
By Proposition~\ref{komponenten}, we have
\[
(\eta_a \cup \eta_b)_{\q_i}= 0 = (\eta_a \cup \eta_b)_{\p_i},
\]
for $a,b, i\in \{1,\ldots, m\} $, hence the cup-product $V\otimes V \to H_\et^2(X\sm S, T)$ is trivial. We claim that the cup-product
\[
H_\et^1(X\sm S_0, T) \otimes V \lang H_\et^2(X\sm S, T)
\]
is surjective. For that purpose we choose elements $\chi_a,\psi_a \in H_\et^1(X\sm S_0, T)$, $a=1,\ldots, m$, such that
\[
\chi_a(\T_{\p_a})\neq 0,\ \chi_a(\Frob_{\q_a})=0,\ \psi_a(\Frob_{\p_a})\neq 0.
\]
We claim that the elements $\chi_1\cup \eta_1,\ldots, \chi_m\cup \eta_m, \psi_1\cup \eta_1,\ldots, \psi_m\cup \eta_m$ generate the $2m$-dimensional vector space $H_\et^2(X\sm S, T) $. In order to see this, we consider the matrix

\smallskip
\[
\left(
\begin{array}{cccccc}
(\chi_1 \cup \eta_1)_{\p_1}& \cdots & (\chi_1 \cup \eta_1)_{\p_m} & (\chi_1 \cup \eta_1)_{\q_1}& \cdots & (\chi_1 \cup \eta_1)_{\q_m}\\
&\vdots && \vdots&&\vdots \\
(\chi_m \cup \eta_m)_{\p_1}&\cdots &(\chi_m \cup \eta_m)_{\p_m}& (\chi_m \cup \eta_m)_{\q_1}& \cdots&(\chi_m \cup \eta_m)_{\q_m}\\
(\psi_1 \cup \eta_1)_{\p_1}& \cdots & (\psi_1 \cup \eta_1)_{\p_m} & (\psi_1 \cup \eta_1)_{\q_1}& \cdots & (\psi_1 \cup \eta_1)_{\q_m}\\
&\vdots && \vdots&&\vdots \\
(\psi_m \cup \eta_m)_{\p_1}&\cdots &(\psi_m \cup \eta_m)_{\p_m}& (\psi_m \cup \eta_m)_{\q_1}& \cdots&(\psi_m \cup \eta_m)_{\q_m}\\
\end{array} \right) \lower1.35cm\hbox{.}
\]

\smallskip\noindent
Denoting a nonzero element by $*$ and an arbitrary element by $?$, Proposition~\ref{komponenten} and our choices show that this matrix has the shape

\smallskip
\[
\left(
\begin{array}{cccccccc}
*&0& \cdots & 0 & 0&0& \cdots & 0\\
0&*& \cdots & 0 & 0&0& \cdots & 0\\
&&\ddots && &&\vdots& \\
0&0&\cdots &*& 0& 0 &\cdots&0\\
?&?& \cdots  &? & *&0& \cdots & 0\\
?&?& \cdots  &? & 0&*& \cdots & 0\\
&&\vdots && \vdots&&\ddots \\
?&?& \cdots  &? & 0&0& \cdots & *\\
\end{array} \right) \lower1.7cm\hbox{.}
\]

\smallskip\noindent
This shows the claim. In particular, the cup-product
\[
H_\et^1(X\sm S, T) \otimes H_\et^1(X\sm S, T) \to H_\et^2(X\sm S, T)
\]
is surjective. As it factors through
\[
H^2(G_S^T(k)(p)) \hookrightarrow H_\et^2(X\sm S, T),
\]
this inclusion is an isomorphism. Therefore the group $G_S^T(k)(p)$ satisfies the assumptions of Theorem~\ref{mildkrit} and we obtain $\cd\, G_S^T(k)(p)=2$. Now Corollary~\ref{kpi1crit} implies the $K(\pi,1)$-property for $(X\sm S, T)$. Finally, by our choice of $S_0$, we have $V_{S\sm \{\p\}}^T(k)=0$ for each $\p\in S$, hence all $\p\in S$ ramify in $k_S^{T,\el}$, in particular, in $k_S^T(p)$.
\end{proof}

For the case $\delta=1$, we need the following  Lemma.

\begin{lemma} \label{einheitswda} Let $\delta=1$ and let $T$ be a finite, nonempty set of places with $V_T^\varnothing(k)=0$ and $S_p\subset T$. We assume $k$ to be totally imaginary if $k$ is a number field (automatic if $p\neq 2$). Then we have  $\null_p \Cl_T(k)=0$ and
\[
k_T^\el = k(\sqrt[p]{E_{k,T}}).
\]
Moreover, let $S=\{\p_1,\ldots,\p_n\}$ be a finite set of places which is disjoint from~$T$. Then we have an inclusion
\[
 k_S^{T,\el} \subset k_T^\el(\sqrt[p]{s_1},\ldots,\sqrt[p]{s_n}),
\]
where $s_i=s_{\p_i}$, $i=1,\ldots,n$. The following assertions are equivalent:

\smallskip
\begin{compactitem}
\item[\rm (i)] the elements $\Frob_{\p_1}, \ldots, \Frob_{\p_n} $ generate $G(k_T^\el|k)$.\smallskip
\item[\rm (ii)]  $V_S^T(k)=0$.
\end{compactitem}

\smallskip\noindent
For a subset $I \subset \{1,\ldots,n\}$ the following are equivalent:

\smallskip
\begin{compactitem}
\item[\rm (a)] The elements $\{\Frob_{\p_i}, i\in I\}$ are linearly independent in $G(k_T^\el|k)$.
\item[\rm (b)] The extensions $k_T^\el$, $k_S^{T,\el}$ and $k(\sqrt[p]{s_i}, i\in I)$ are linearly disjoint over $k$.
\end{compactitem}
\end{lemma}

\begin{remark}
If $\Frob_{\p_1}, \ldots, \Frob_{\p_n} $ are linearly independent in  $G(k_T^\el|k)$, then the inclusion $k_S^{T,\el} \subset k_T^\el(\sqrt[p]{s_1},\ldots,\sqrt[p]{s_n})$ seems to contradict assertion (b) for $I=\{1,\ldots,n\}$. But we have $k_S^{T,\el}=k$ in this case.
\end{remark}

\begin{proof}[Proof of Lemma~\ref{einheitswda}] By Theorem~\ref{bstdual}, we have $\Sha^2(k,T)=0$. Hence Poitou-Tate duality (\cite{NSW}, Theorem 8.6.7) implies $\Sha^1(k,T,\mu_p)=0$. Since $\delta=1$ and by \cite{NSW}, Lemma~8.6.3, we obtain $\Hom(\Cl_T(k),\Z/p\Z)= \Sha^1(k,T)=0$. As $\Cl_T(k)$ is finite, this implies $\null_p\Cl_T(k)=0$. By Dirichlet's unit theorem, we have  $\dim_{\F_p} E_{k,T}/p = \# T + r$  and  Theorem~\ref{globcoh} implies
\[
h^1(X\sm T)= \# T + \sum_{\p\in S_p} [k_\p:\Q_p] - r = \# T +r.
\]
This shows that $k(\sqrt[p]{E_{k,T}})=k_T^\el$.
Lemma~\ref{VSchange} yields an isomorphism
\[
E_{k,T}/p \liso V_\varnothing^T(k).
\]
Now let $S$ be a finite set of places which is disjoint from  $T$ and let $k(\sqrt[p]{\alpha})$, $\alpha \in k^\times/k^{\times p}$, be a cyclic subextension of $k_S^{T,\el}|k$. Then we have $\alpha\in V_T^S(k)$ and we find exponents $a_1,\ldots,a_n$ such that
$\alpha \cdot \bar s_1^{a_1}\cdots \bar s_n^{a_n} \in V_\varnothing^T(k)= E_{k,T}/p$. Hence $k(\sqrt[p]{\alpha}) \subset k(\sqrt[p]{E_{k,T}},\sqrt[p]{s_1},\ldots,\sqrt[p]{s_n})$ and therefore
\[
k_S^{T,\el} \subset k_T^\el(\sqrt[p]{s_1},\ldots,\sqrt[p]{s_n}).
\]
Next we prove the equivalence of (i) and (ii). We have $V_S^T(k)=0$ if and only if the natural map  $E_{k,T}/p \to \prod_{i=1}^n k_{\p_i}^\times/k_{\p_i}^{\times p}$ is injective. This is true if and only if for no element $e\in E_{k,T}/p$, $e\neq 1$, the cyclic extension $k(\sqrt[p]{e})|k$ is completely split at~$S$. Hence (ii) holds if and only if the elements $\Frob_{\p_1}, \ldots, \Frob_{\p_n}$ generate the Galois group $G(k(\sqrt[p]{E_{k,T}})|k)=G(k_T^\el|k)$.  Therefore  (i) and (ii) are equivalent.

\medskip
We denote by $I_{k,T}$ the group of  $T$-id\`{e}les and by $C_T(k)$ the group of $T$-id\`{e}le classes of~$k$. Since  $\Cl_T(k)(p)=0$, \cite{NSW}, Proposition 8.3.5, implies the exact sequence
\[
0 \lang E_{k,T} \otimes \Z_p \lang I_{k,T} \otimes \Z_p \lang C_T(k) \otimes \Z_p \lang 0\,.
\]
Now we consider a subset  $I\subset \{1,\ldots, n\}$. Let  $H_I \subset k^\times/k^{\times p}$ be the subgroup generated by the elements $s_i$, $i\in I$. Since $S_p\subset T$, we have
\[
k_S^{T,\el}=k\left(\sqrt[p]{V_T^S(k)}\right).
\]
By Kummer theory, the extensions $k_T^\el$, $k_S^{T,\el}$ and $k(\sqrt[p]{s_i}, i\in I)$ are linearly disjoint over $k$ if and only if the homomorphism
\[
E_{k,T}/p \times V_T^S(k) \times H_I  \lang  k^\times/k^{\times p}
\]
is injective. Because of $H_I \cap E_{k,T}/p = 1$, this is equivalent to
\[
( H_I \cdot E_{k,T}/p ) \cap V_T^S(k) =  1
\]
and hence to the injectivity of the map
\[
H_I \cdot E_{k,T}/p \lang \prod_{\p\in T} k_\p^\times/k_\p^{\times p}=I_{k,T}/p.
\]
This latter map is injective if and only if the composite map
\[
H_I \lang (I_{k,T}/p)/ \im (E_{k,T}/p) \liso C_T(k)/p \mathop{\longrightarrow}\limits_\rec^\sim G(k_T^\el|k)
\]
is injective. By class field theory, $s_i$ maps to $\Frob_{\p_i}^{-1}\in G(k_T^\el|k)$. Therefore the injectivity of the composite map above is equivalent to assertion~(a).
\end{proof}

\begin{proof}[Proof of Theorem~\ref{hilf} in the case $\delta=1$] We choose a sufficiently large finite set of places $T_0$ such that $V_{T\cup T_0}^\varnothing(k)=0$ and $S_p\subset T\cup T_0$. To simplify the notation, we replace $T$ by $T\cup T_0$ in what follows.
Now we choose to each nonzero element $g\in G(k_T^\el|k)$ a prime $\p_g$  with $\p\notin T\cup \M$ and $g=\Frob_{\p_g}$. We denote the set of these primes by ${S_0}$ and we choose a numbering
\[
{S_0}=\{\p_0, \ldots, \p_m\}.
\]
As before, we put $s_i=s_{\p_i}$. By Lemma~\ref{einheitswda}, we have  $V_{S_0}^T(k)=0$, $k_T^\el= k(\sqrt[p]{E_{k,T}})$ and $ k_{S_0}^{T,\el} \subset k_T^\el(\sqrt[p]{s_0},\ldots,\sqrt[p]{s_m})$. By Theorem~\ref{globcoh}, we obtain $h^2(X\sm {S_0},T)=m$ and $h^1(X\sm T)=\# T +r:=n$.

\smallskip
Let $a$, $1\leq a\leq m$, be an index. We choose a subset $I_a \subset \{1,\ldots,m\}$ of cardinality $n-1$ with $a\notin I_a$ such that
$(\Frob_{\p_0}, \{ \Frob_{\p_i}\}_{i\in I_a})$ as well as $(\Frob_{\p_a}, \{ \Frob_{\p_i}\}_{ i\in I_a})$ is a basis of $G(k_T^\el|k)$. This possible. Indeed, if $\Frob_{\p_0}$ and $\Frob_{\p_a}$ are linearly dependent in $G(k_T^\el|k)$, then we complete $\Frob_{\p_0}$ to a basis. If these two elements are linearly independent, then we choose  $a'$ with $\Frob_{\p_{a'}}=\Frob_{\p_0} + \Frob_{\p_a}$ and complete $(\Frob_{\p_0}, \Frob_{\p_{a'}})$ to a basis.

By Lemma~\ref{einheitswda}, the extensions  $k_T^\el|k$ (degree $p^n$), $k_{S_0}^{T,\el}$ (degree $p^{m+1-n}$) and $k(\sqrt[p]{s_a}, \sqrt[p]{s_i}, i\in I_a)$ (degree $p^n$) are linearly disjoint over $k$. By reasons of degree, their composite is equal to  $k_T^\el(\sqrt[p]{s_0},\ldots,\sqrt[p]{s_m})$ (degree $p^{m+1+n}$).

By our choice of $I_a$, the same is true after replacing  $k(\sqrt[p]{s_a}, \sqrt[p]{s_i}, i\in I_a)$ by $k(\sqrt[p]{s_0}, \sqrt[p]{s_i}, i\in I_a)$.

For $i\in \{0,\ldots,m\}$,  let $\T_i \subset G(k_{S_0}^{T,\el}|k)$ be the inertia group of~$\p_i$. Since the elements  $\Frob_{\p_j}$, $j\neq i$, still generate $G(k_T^\el|k)$, we have $V_{{S_0}\backslash \{\p_i\}}^T(k)=0$ by Lemma~\ref{einheitswda}. Therefore $\p_i$ ramifies in $k_{S_0}^{T,\el}|k$ and $\T_i$ is cyclic of order~$p$. By construction, the  $m-n+1$ cyclic subgroups of order~$p$
\[
\T_i, \ i \notin \{ a\} \cup I_a\,,
\]
generate the $(m-n+1)$-dimensional vector space $G(k_{S_0}^{T,\el}|k)$ and the same is true for the subgroups
\[
\T_i, \ i \notin \{ 0\} \cup I_a\,.
\]
Therefore the subgroups $\T_i$,  $i \notin \{ 0, a\} \cup I_a$ generate an $(m-n)$-dimensional subspace, the extension $k_{\{\p_0,\p_a\}}^{T,\el}$ is cyclic of order~$p$ and ramified at $\p_0$ and $\p_a$.

\noindent
For $a=1,\ldots, m$ and a prime $\q$, we consider the following condition $(C_a)$:

\medskip
\begin{compactitem}
\item $\q \notin T\cup \M$.\smallskip
\item $\q$ splits completely in $k_T^\el|k$.\smallskip
\item For each $i\in I_a$,  $\q$ splits in  $k(\sqrt[p]{s_i})|k$.\smallskip
\item $\q$ is inert in the extension $k(\sqrt[p]{s_a})|k$.\smallskip
\item the image of  $\Frob_{\q}$ in $G(k_{S_0}^{T,\el}|k)$ is  in $G(k_{S_0}^{T,\el}|k_{\{\p_0,\p_a\}}^{T,\el})\sm \{ 0\}$.\smallskip
\end{compactitem}

\medskip\noindent
Applying the Chebotarev density theorem to the elementary abelian extension
\[
k_T^\el(\sqrt[p]{s_0},\ldots,\sqrt[p]{s_m} )|k,
\]
we find a prime $\q_1$ in $k$ which satisfies condition $(C_1)$.  By  Lemma~\ref{dual}, the extension $k_{\{\q_1\}}^{T,\el}$ is cyclic of degree $p$ and ramified at $\q_1$. Now we choose recursively primes $\q_2,\ldots,\q_m$ in $k$ such that

\medskip
\begin{compactitem}
\item $\q_a$ satisfies condition $(C_a)$, and \smallskip
\item For $b<a$, $\q_a$ splits in $k_{\{\q_b\}}^{T,\el}|k$.
\end{compactitem}

\medskip\noindent
In particular, the $\q_i$ are pairwise different and we have $N(\q_i)\equiv 1 \bmod p$. We claim that
\[
S=\{\p_0,\ldots,\p_m,\q_1,\ldots,\q_m\}
\]
has the required properties. We have  $h^2(X\sm S,T)=2m$ and, by Corollary~\ref{h2lokali}, the natural map
\[
H_\et^2(X\sm S, T) \lang \prod_{i=1}^m H^2(k_{\p_i}) \oplus \prod_{i=1}^m H^2(k_{\q_i} )
\]
is an isomorphism.

\noindent
Let $\eta_a$ be a generator of $H_\et^1(X\sm \{ \q_a\}, T)$. We consider the $m$-dimensional vector space $V$ spanned by $\eta_1,\ldots,\eta_m$ in $H_\et^1(X\sm S, T)$. Then we have
\[
H_\et^1(X\sm S, T) \cong H_\et^1(X\sm S_0, T) \oplus V.
\]
Lemma~\ref{dual} implies
\[
k_T^\el k_{\{\q_a\}}^{T,\el} = k_T^\el(\sqrt[p]{s_{\q_a}}).
\]
By Proposition~\ref{komponenten}, we have
\[
(\eta_a \cup \eta_b)_{\q_i}= 0 = (\eta_a \cup \eta_b)_{\p_i}
\]
for $a,b, i\in \{1,\ldots, m\} $. It remains to show that the cup-product
\[
H_\et^1(X\sm S_0, T) \otimes V \lang H_\et^2(X\sm S, T)
\]
is surjective. To this end we choose for each $a\in \{1,\ldots, m\}$ a generator $\chi_a$ of $H_\et^1(X\sm \{\p_0,\p_a\}, T)$. By construction, we have  $\chi_a(\T_a)\neq 0$, $\chi_a(\T_i)=0$ for $i\notin \{0,a\}\cup I_a$ and $\chi_a(\Frob_{\q_a})=0$. Furthermore, we choose $\psi_a\in H_\et^1(X\sm S, T)$ with $\psi_a(\Frob_{\q_a})\neq 0$.

We claim that the elements  $\chi_1\cup \eta_1,\ldots, \chi_m\cup \eta_m, \psi_1\cup \eta_1,\ldots, \psi_m\cup \eta_m$ generate the $2m$-dimensional vector space $H_\et^2(X\sm S, T) $. This can be seen in exactly the same way as in the case $\delta=0$ and also the rest of the proof is word by word the same from this point on.
\end{proof}

\section{Proof of Theorem~\ref{haupt}}

In this section we prove Theorem~\ref{haupt}.
We introduce the following notation: let  $K|k$ be a (typically infinite) separable algebraic extension and let $S$ be a finite set of places of $k$. Then we write (as before omitting the coefficients $\F_p$)
\[
\ressum_{{\mathfrak p} \in S(K)} H^i(K_\p)
\stackrel{df}{=}
\varinjlim_{k'\subset K} \bigoplus_{{\mathfrak p} \in S(k')} H^i(k'_\p),
\]
where the limit on the right hand side runs through all finite subextensions $k'|k$ of $K|k$.
If $K|k$ is Galois with group  $G(K|k)$, then the generalized sum $\ressumsmall_{{\mathfrak p} \in S(K)} H^i(K_\p)$ is the maximal discrete $\Gal(K|k)$-submodule of the product  $\prod_{{\mathfrak p} \in S(K)} H^i(K_\p)$.

\medskip
The next proposition reduces Theorem~\ref{haupt} to Theorem~\ref{hilf}.

\begin{proposition}\label{erweiterung}
Let $S$, $T$, $S_0$ and $T_0$ be pairwise disjoint finite sets of places of the global field~$k$. Assume that $(X\sm S_0, S\cup T \cup T_0)$ has the $K(\pi,1)$-property for $p$. Then also $(X\sm S_0, S\cup T)$, $(X\sm S_0, T)$ and $(X\sm (S\cup S_0), T)$ have the  $K(\pi,1)$-property for $p$. Furthermore, the extension $k_{S\cup S_0}^T(p)$ realizes the maximal $p$-extension $k_\p(p)$ for all $\p\in S$ and the natural homomorphism
\[
\freeproductmed_{\p\in S(k_{S_0}^{S\cup T}(p))}G(k_\p(p)|k_\p) \lang G\big(k_{S\cup S_0}^T(p)|k_{S_0}^{S\cup T}(p)\big).
\]
is an isomorphism.
\end{proposition}

\begin{proof} To simplify notation, we write $k_S^T$ for $k_S^T(p)$. We start by calculating the cohomology of  $(X\sm S_0, S\cup T)_{k_{S_0}^{S\cup T\cup T_0}}$.
Since $(X\sm S_0, S\cup T \cup T_0)$ has the $K(\pi,1)$-property for~$p$,  Lemma~\ref{kpi1lem} implies
\[
H^i_\et\big((X\sm S_0, S\cup T \cup T_0)_{k_{S_0}^{S\cup T \cup T_0}}\big)=0,\quad i\geq 1.
\]
Therefore Corollary~\ref{ausschneid2} shows that $H^i_\et((X\sm S_0,S\cup T)_{k_{S_0}^{S\cup T \cup T_0}})=0$ for $i\geq 2$ and we obtain an isomorphism
\[
H^1_\et\big((X\sm S_0,S\cup T)_{k_{S_0}^{S\cup T \cup T_0}}\big)\liso
 \ressum_{\p\in T_0(k_{S_0}^{S\cup T \cup T_0})} \back
H^1_\nr(k_\p).
\]
In particular, $k_{S_0}^{S\cup T}$ realizes the maximal elementary abelian unramified $p$-extension $k_\p^{\nrel}$ of $k_\p$ for all $\p\in T_0$. The Hochschild-Serre \pagebreak spectral sequence implies an inclusion
\[
H^2\big(G(k_{S_0}^{S\cup T}|k_{S_0}^{S\cup T\cup T_0})\big)\hookrightarrow H^2_\et\big((X\sm S_0,S\cup T)_{k_{S_0}^{S\cup T \cup T_0}}\big)=0,
\]
showing that the pro-$p$-group $G(k_{S_0}^{S\cup T}|k_{S_0}^{S\cup T\cup T_0})$ is free. Therefore Lemma~\ref{kpi1crit} shows that $(X\sm S_0,S\cup T)_{k_{S_0}^{S\cup T \cup T_0}}$ has the $K(\pi,1)$-property for~$p$. We obtain
\[
H^i_\et\big((X\sm S_0,S\cup T)_{k_{S_0}^{S\cup T}}\big)=0,\quad i\geq 1,
\]
showing the $K(\pi,1)$-property for  $(X\sm S_0, S\cup T)$ by Lemma~\ref{kpi1lem}. The decomposition groups  $Z_\p(k_{S_0}^{S\cup T}|k_{S_0}^{S\cup T\cup T_0})$ of the (unramified) places $\p\in T_0$ are nontrivial and torsion-free. Hence the field $k_{S_0}^{S\cup T}$ realizes the maximal unramified $p$-extension $k_\p^\nr(p)$ of $k_\p$ for all $\p\in T_0$.

The same arguments also show that $(X\sm S_0, T)$ has the  $K(\pi,1)$-property for~$p$ and that $k_{S_0}^T$ realizes the maximal unramified $p$-extension $k_\p^\nr(p)$ of $k_\p$ for all $\p\in S$. Using the excision sequence, we obtain isomorphisms \[
H^i_\et\big((X\sm (S\cup S_0),T)_{k_{S_0}^{T}}\big)\liso
\ressum_{\p\in S(k_{S_0}^{ T})}  H^i_{/\nr}(k_\p),\quad i\geq 1.
\]
As $k_{S_0}^{T}$ realizes the maximal unramified $p$-extensions $k_\p^\nr(p)$ of the places $\p\in S$, these cohomology groups vanish for $i\geq 2$ and $k_{S\cup S_0}^T$ realizes the maximal elementary abelian $p$-extension of~$k_\p^\nr(p)$ for all  $\p\in S$. As above, we obtain that the curve $(X\sm(S\cup S_0), T)$ has the $K(\pi,1)$-property for $p$ and that the pro-$p$-group $G(k_{S\cup S_0}^T|k_{S_0}^T)$ is free.  Therefore the natural homomorphisms
\[
G\big(k_\p(p)|k_\p^\nr(p)\big) \longrightarrow Z_\p(k_{S\cup S_0}^T|k_{S_0}^T),\quad \p\in S(k_{S_0}^T),
\]
from the full local groups onto the decomposition groups are homomorphisms between free pro-$p$-groups which induce isomorphisms on $H^1(-,\F_p)$.  Hence they are isomorphisms (see \cite{NSW}, 1.6.15). We conclude that $k_{S\cup S_0}^T$ realizes the maximal $p$-extension $k_\p(p)$ for all $\p\in S$. Now we consider another excision sequence to obtain isomorphisms
\[
H^i_\et\big((X\sm (S\cup S_0),T)_{k_{S_0}^{S\cup T}}\big)\liso \ressum_{\p\in S(k_{S_0}^{S\cup T})} \back H^i(k_\p).
\]
Since $(X\sm (S\cup S_0), T)$ has the $K(\pi,1)$-property for~$p$, the cohomology of the group $G(k_{S\cup S_0}^T|k_{S_0}^{S\cup T})$ coincides with the \'{e}tale cohomology of  $(X\sm (S\cup S_0),T)_{k_{S_0}^{S\cup T}}$. Using the calculation of the cohomology of a free product (\cite{NSW}, Theorem 4.3.14), we conclude that
\[
\phi: \freeproductmed_{\p\in S(k_{S_0}^{S\cup T})}G(k_\p(p)|k_\p) \lang G(k_{S\cup S_0}^T|k_{S_0}^{S\cup T})
\]
is a homomorphism between pro-$p$-groups which induces isomorphisms on $H^i(-,\F_p)$ for all~$i$. By \cite{NSW}, Proposition 1.6.15, $\phi$ is an isomorphism.
\end{proof}

Now we deduce Theorem~\ref{haupt}. Let $S$, $T$ and $\M$ be pairwise disjoint sets of places of the global field $k$, where $S$ and $T$ are finite and $\M$ has Dirichlet density $\delta(\M)=0$.
We choose sets of places $S_0$, $T_0$ to $S\cup T$ and $\M$ as in Theorem~\ref{hilf}, i.e.\

\medskip
\begin{compactitem}
\item $S_0$ is a nonempty set of places $\p$ of norm $N(\p)\equiv 1\bmod p$, \smallskip
\item $S_0\cap (S\cup T\cup T_0\cup \M)=\varnothing$, \smallskip\pagebreak
\item $(X\sm S_0, S\cup T\cup T_0)$ has the $K(\pi,1)$-property for~$p$, \smallskip
\item each $\p\in S_0$ ramifies in $k_{S_0}^{S\cup T\cup T_0}(p)$, \smallskip
\item $V_{S_0}^{S\cup T\cup T_0}(k)=0$,\smallskip
\item the cup-product
\[
H^1\big(G_{S_0}^{S\cup T\cup T_0}(k)(p)\big) \otimes H^1\big(G_{S_0}^{S\cup T\cup T_0}(k)(p)\big) \lang H^2\big(G_{S_0}^{S\cup T\cup T_0}(k)(p)\big)
\]
is surjective.
\end{compactitem}

\medskip\noindent
By Proposition~\ref{erweiterung}, we obtain the following assertions of Theorem~\ref{haupt}:
we have $cd\, G_{S\cup S_0}^T(p) \leq 2$, assertion (ii) but at the moment only for the places in $S$, and assertions (iii) and (iv).

Now we show that  $k_{S\cup S_0}^T(p)$ realizes the  maximal $p$-extension $k_\p(p)$ of $k_\p$ for all places $\p\in S_0$. Let  $\p\in S_0$. Then $\p$ is not above $p$ and the local field $k_\p$ contains a primitive $p$-th root of unity. The decomposition group $Z_\p(k_{S\cup S_0}^T(p)|k)$ has cohomological dimension less or equal to~$2$ as it is a subgroup of $G(k_{S\cup S_0}^T(p)|k)$. Passing in \cite{NSW}, Theorem~7.5.2, to the maximal pro-$p$-factor group, we see that the full local group  $\Gal(k_\p(p)|k_\p)$ is, as a pro-$p$-group, generated two generators $\sigma, \tau$ with the one defining relation $\sigma\tau\sigma^{-1}=\tau^q$. The element $\tau$ is a generator of the inertia group,  $\sigma$ is a Frobenius lift and  $q=N(\p)$. Therefore $\Gal(k_\p(p)|k_\p)$ has exactly three factor groups of cohomological dimension less or equal to $2$: itself, the trivial group and the Galois group of the maximal unramified $p$-extension
of~$k_\p$. Since $\p$ ramifies in the extension  $k_{S\cup S_0}^T(p)$, the decomposition group must be full, i.e.\ $k_{S\cup S_0}^T(p)|k$ realizes the maximal $p$-extension $k_\p(p)|k_\p$. Finally, $S_0$ is nonempty and for $\p\in S_0$ the decomposition group $Z_\p(k_{S\cup S_0}^T(p)|k)$ has cohomological dimension~$2$. Hence also $G(k_{S\cup S_0}^T(p)|k)$ has cohomological dimension~$2$.
Since $V_{S_0}^{S\cup T\cup T_0}(k)=0$, also the groups $V_{S_0}^{S\cup T}(k)$ and $V_{S\cup S_0}^{T}(k)$ vanish. Therefore the sequences
\[
0\to H_\et^2(X\sm S_0, S\cup T) \to H_\et^2(X\sm (S\cup S_0),  T)\to \bigoplus_{\p\in S} H_\et^2(k_\p)\to 0,
\]
\[
0\to H_\et^1(X\sm S_0, S\cup T) \to H_\et^1(X\sm (S\cup S_0),  T)\to \bigoplus_{\p\in S} H_\et^1(k_\p)\to 0
\]
are exact. Furthermore, the cup-product
\[
H_\et^1(X\sm S_0, S\cup T) \otimes H_\et^1(X\sm S_0, S\cup T) \lang H_\et^2(X\sm S_0, S\cup T)
\]
is surjective. This follows from our choice of  $S_0$ and since the homomorphism
\[
 H_\et^2(X\sm S_0, S\cup T\cup T_0) \lang H_\et^2(X\sm S_0, S\cup T).
\]
is surjective. Finally, recall that the local cup-products $H_\et^1(k_\p)\otimes H_\et^1(k_\p)\to H_\et^2(k_\p)$ are  surjective and that the inflation homomorphisms $H^i(G(k_\p(p)|k_\p))\to H^i(k_\p)$ are isomorphisms for all~$i$ (see \cite{NSW}, 7.5.8). As $k_{S\cup S_0}^T(p)$ realizes the maximal $p$-extension $k_\p(p)$ of $k_\p$ for $\p\in S$, we obtain the surjectivity of the cup-product
\[
H_\et^1(X\sm (S\cup S_0),  T) \otimes H_\et^1(X\sm (S\cup S_0), T) \lang H_\et^2(X\sm (S\cup S_0), T)
\]
This concludes the proof of  Theorem~\ref{haupt}.

\pagebreak

\section{Enlarging the set of places} \label{erwsec}

\begin{theorem} \label{enlarge} Let $T$ and $S'$ be disjoint, finite sets of places of the global field~$k$, $S\subset S'$  a subset and $p\neq \text{char}(k)$ a prime number. Suppose $p\neq 2$ or $k$ totally imaginary if $k$ is a number field. Assume that the marked arithmetic curve $(X\sm S,T)$ has the $K(\pi,1)$-property for~$p$. If no prime  $\p\in S'\sm S$ splits completely in the extension  $k_S^T(p)|k$, then the following hold:

\medskip
\begin{compactitem}
\item[{\rm (i)}] Also $(X\sm S',T)$ has the $K(\pi,1)$-property for $p$.\smallskip
\item[{\rm (ii)}] $k_{S'}^T(p)_\p=k_\p(p)$ for all $\p \in S'\sm S$.
\end{compactitem}

\medskip\noindent
Furthermore, the arithmetic form of Riemann's existence theorem holds, i.e.\ for \linebreak $K=k_S^T(p)$, the natural homomorphism \[
\freeproductmed_{{\mathfrak p} \in S'\backslash S(K)} \Gal(K_\p(p)|K_\p) \longrightarrow \Gal(k_{S'}^T(p)|K)
\]
is an isomorphism. In particular,  $\Gal(k_{S'}^T(p)|k_S^T(p))$ is a free pro-$p$-group.
\end{theorem}

\begin{remark}
If $k_S^T(p)|k$ is infinite, then the set of places which split completely has Dirichlet density~zero. This statement can be sharpened in the style of \cite{Ih}, see \cite{TV}, Proposition~3.1. Naturally the question arises, whether this set is finite or even equal to~$T$. This is clear in the case $S\supset S_p$, $T=\varnothing$, because then the cyclotomic  $\Z_p$-extension of~$k$ is contained in~$k_S^T(p)$.
\end{remark}

\begin{proof}[Proof of Theorem~\ref{enlarge}] The $K(\pi,1)$-property implies
\[
H^i(G_S^T(k)(p),\F_p) \cong H^i_\et(X\sm S,T,\F_p)=0 \text{ for } i\geq 4,
\]
in particular, we have $\cd\, G_S^T(k)(p)\leq 3$. Let $\p \in S'\sm S$ and put $K=k_S^T(p)$. By assumption, $\p$ does not split completely in~$K|k$.  Because of  $\cd\, G_S^T(k)(p)< \infty$, the decomposition group of~$\p$ in $K|k$ is a nontrivial and torsion-free factor group of $\Z_p\cong \Gal(k_\p^{\nr}(p)|k_\p)$. Therefore
\[
K_\p=k_\p^{\nr}(p)
\]
for all $\p\in S'\sm S$. We consider the excision sequence for the curves $(X \sm S,T)_{K}$ and $(X \sm S',T)_{K}$. Since  $(X\sm S,T)$ has the $K(\pi,1)$-property for~$p$, we obtain $H^i_{\et}((X \sm S,T)_{K},\F_p)=0$ for $i \geq 1$. Omitting the coefficients $\F_p$ from the notation, we obtain isomorphisms
\[
H^i_{\et}\big((X \sm S',T)_K\big) \stackrel{\sim}{\longrightarrow} \ressum_{{\mathfrak p} \in S'\backslash S (K)} H^{i+1}_{\mathfrak p}\big((X \sm S,T)_{K}\big)
\]
for $i\geq 1$. This implies
\[
H^i_{\et}\big((X\sm S',T)_{K}\big)=0\
\]
for $i\geq 2$. Now $ (X\sm S',T)_{k_{S'}^T(p)}$ is the universal pro-$p$-covering of the marked curve
$(X\sm S',T)_{K}$. Therefore the Hochschild-Serre spectral sequence
\[
E_2^{ij}=H^i\big(\Gal(k_{S'}^T(p)|K), H^j_\et((X\sm S',T)_{k_{S'}^T(p)})\big) \Rightarrow H^{i+j}_\et((X\sm S',T)_{K})
\]
yields an inclusion
\[
H^2\big(\Gal(k_{S'}^T(p)|K)\big) \hookrightarrow H^2_{\et}\big((X\sm S',T)_{K}\big)=0.
\]
Hence  $\Gal(k_{S'}^T(p)|K)$ is a free pro-$p$-group and we have an isomorphism
\[
H^1(\Gal(k_{S'}^T(p)|K)) \stackrel{\sim}{\to} H^1_{\et}\big((X\sm S',T)_{K}\big)
\cong \ressum_{{\mathfrak p}\in S'\backslash S(K)} H^1(K_{\mathfrak p}).
\]
We consider the natural homomorphism
\[
\phi: \freeproductmed_{{\mathfrak p} \in S'\backslash S(K)} \Gal(K_\p(p)|K_\p) \longrightarrow \Gal(k_{S'}^T(p)|K).
\]
Because of $K_\p=k_\p^\nr(p)$ for $\p\in S'\sm S$, the factors in the free product on the left hand side are free pro-$p$-groups. By the calculation of the cohomology of a free product (\cite{NSW}, 4.3.10 and 4.1.4), $\phi$ is a homomorphism between free pro-$p$-groups which induces an isomorphism on $H^1(-,\F_p)$.  Therefore $\phi$ is an isomorphism (see \cite{NSW}, 1.6.15).  In particular, we have $k_{S'}^T(p)_\p=k_\p(p)$ for all $\p\in S'\sm S$.
Using the freeness of  $\Gal(k_{S'}^T(p)|K)$, the Hochschild-Serre spectral sequence induces an isomorphism
\[
0=H^2_\et((X\sm S',T)_{K})\mapr{\sim} H^2_\et\big((X\sm S',T)_{k_{S'}(p)}\big)^{\Gal(k_{S'}^T(p)|K)}.
\]
Since $\Gal(k_{S'}(p)|k_S(p))$ is a pro-$p$-group, this implies  $H^2_\et\big((X\sm S',T)_{k_{S'}(p)}\big)=0$. Therefore $(X\sm S',T)$ has the $K(\pi,1)$-property for~$p$ by Lemma~\ref{kpi1lem}.
\end{proof}

\section{Duality for the fundamental group}\label{dualsec}
We start by investigating the connection between the $K(\pi,1)$-property and the universal norms of global units.

First, we remove redundant places from $S$: if $\p \nmid p$ is a place with $\zeta_p\notin k_\p$, then every  $p$-extension of the local field~$k_\p$ is unramified (see \cite{NSW}, 7.5.9).  Therefore places  $\p \notin  S_p$ with $N(\p) \not \equiv 1 \bmod p$ cannot ramify in a $p$-extension.
Removing all these redundant primes from $S$, we obtain a subset  $S_{\min} \subset S$ with $G_S^T(p)=G_{S_{\min}}^T(p)$. If $\delta=1$, then  $S_\min=S$, i.e.\ there are no redundant places.

\begin{lemma}
$(X\sm S,T)$ has the $K(\pi,1)$-property for $p$ if and only if this is the case for  $(X\sm S_\min,T)$.
\end{lemma}

\begin{proof}
We have $k_S^T(p)=k_{S_\min}^T(p)$. Denoting this field by $K$, Proposition~\ref{localcoh} implies
\[
H^i_\p\big((X,T)_K,\F_p\big)=0
\]
for $i\geq 1$ and every $\p \in S\sm S_\min(K)$. The excision sequence shows that, for $i\geq 1$, the group
$H^i_\et((X\sm S,T)_K,\F_p)$ vanishes if and only if  $H^i_\et((X\sm S_\min,T)_K,\F_p)$ vanishes. Therefore the statement of the lemma follows from Lemma~\ref{kpi1lem}.
\end{proof}

\begin{proposition} \label{thmb} Let $S$ and $T$ be disjoint, finite sets of places of the global field~$k$ and let $p\neq \text{char}(k)$ be a prime number. Suppose $p\neq 2$ or $k$ totally imaginary if $k$ is a number field.  Then any two of the conditions {\rm (a)--(c)} below imply the third.

\medskip\noindent
\begin{compactitem}
\item[\rm (a)] $(X\sm S,T)$ has the $K(\pi,1)$-property for $p$. \smallskip
\item[\rm (b)] $\varprojlim_{K\subset k_S^T(p)} E_{K,T} \otimes \Z_p=0$. \smallskip
\item[\rm (c)] $k_{S}^T(p)_\p=k_\p(p)$ for all $\p\in S_\min$.
\end{compactitem}

\medskip\noindent
The limit in {\rm (b)} runs through all finite subextensions $K|k$ of $k_S^T(p)|k$.
If {\rm (a)--(c)} hold, then we also have
\[
\varprojlim_{K\subset k_S^T(p)} E_{K,S_\min\cup T} \otimes \Z_p=0.
\]
\end{proposition}

\begin{remarks} 1. By Theorem~\ref{haupt}, (a)--(c) hold after adding finitely many primes to~$S$.

\noindent
2. Assume $\zeta_p \in k$, $S_p\subset S$ and $T=\varnothing$. Then (a) holds and condition (c) holds for  $p>2$ if $\# S > r_2+2$ (see \cite{NSW} 10.9.1 and Remark~2 after 10.9.3). In the case  $k=\Q(\zeta_p)$, $S=S_p$, $T=\varnothing$, condition (c) holds if and only if $p$ is an irregular prime number.
\end{remarks}

\begin{proof}[Proof of Proposition~\ref{thmb}] Without restriction, we may assume $S=S_\min$. Applying the topological  Nakayama lemma (\cite{NSW}, 5.2.18) to the compact $\Z_p$-module  $\varprojlim E_{K,T} \otimes \Z_p$, we see that condition (b) is equivalent to condition (b') below:

\medskip
\begin{compactitem}
\item[\rm (b)']\quad  $\varprojlim_{K\subset k_S^T(p)} E_{K,T}/p=0$.
\end{compactitem}

\medskip\noindent
Furthermore, by Lemma~\ref{kpi1lem}, condition (a) is equivalent to

\medskip
\begin{compactitem}
\item[\rm (a)'] \quad $\varinjlim_{K\subset k_S^T(p)} H^i_\et((X\sm S, T)_K, \F_p)=0$ for $i\geq 1$.
\end{compactitem}

\medskip\noindent
By Theorem~\ref{globcoh},  (a)' holds for $i=1$, $i\geq 4$, and also for  $i=3$ if $S$ is nonempty or  $\delta=0$. If   $S=\varnothing$ and $\delta=1$, then we have $H^3_\et((X\sm S,T)_K,\F_p)\cong \mu_p^\vee$. Therefore, in this case, condition  (a)' for $i=3$ is satisfied if and only if the group $G_S^T(k)(p)$ is infinite.

For every  $K\subset k_S^T(p)$, Lemma~\ref{VSchange} implies the exact sequence
\[
0\longrightarrow E_{K,T} /p \longrightarrow V_\varnothing^T(K) \longrightarrow \null_p \Cl_T(K) \longrightarrow 0\,.
\]
The field $k_S^T(p)$ does not admit unramified $p$-extensions in which all primes of $T$ split completely. Hence
class field theory implies
\[
\varprojlim_{K\subset k_S^T(p)} \null_p\Cl_T(K) \subset \varprojlim_{K\subset k_S^T(p)} \Cl_T(K)\otimes \Z_p=0.
\]
Therefore Theorem~\ref{bstdual} provides an isomorphism
\[
\varinjlim_{K\subset k_S^T(p)} H^2_\et\big((X, T)_K, \F_p\big)=\varinjlim_{K\subset k_S^T(p)} \Sha^2(K,\varnothing, T) \cong \big(\varprojlim_{K\subset k_S^T(p)} E_{K,T}/p \ \big)^\vee. \leqno (*)
\]
We start by showing that conditions (a) and (b) are equivalent in the case  $S=\varnothing$ (where (c) holds for trivial reasons). If (a)' holds, then  $(\ast)$ implies (b)'. If (b) holds, then, in particular,  $\delta=0$ or $G_S^T(k)(p)$ is infinite. Hence we obtain (a)' for $i=3$. Furthermore, (a)' for $i=2$ follows from $(\ast)$ and (b)'. This finishes the proof in the case  $S=\varnothing$.

\medskip
From now on let $S\neq \varnothing$. For  $\p\in S=S_\min$, every proper Galois subextension of  $k_\p(p)|k_\p$ admits ramified $p$-extensions. By the calculation of the local cohomology groups in Proposition~\ref{localcoh},  (c) is therefore equivalent to

\medskip
\begin{compactitem}
\item[\rm (c)'] \quad $\varinjlim_{K\subset k_S^T(p)} \bigoplus_{\p \in S(K)} H^i_\p((X,T)_K, \F_p)=0$ for all $i$,
\end{compactitem}

\medskip\noindent
and also to
\medskip
\begin{compactitem}
\item[\rm (c)''] \quad $\varinjlim_{K\subset k_S^T(p)} \bigoplus_{\p \in S(K)} H^2_\p((X,T)_K, \F_p)=0$.
\end{compactitem}

\medskip\noindent
Now we consider the direct limit over $K\subset k_S^T(p)$ of the excision sequences (coefficients $\F_p)$
\[
\cdots \to  \bigoplus_{\p\in S(K)} H^i_\p\big((X,T)_K\big) \to H^i_\et\big((X,T)_K\big) \to H^i_\et\big((X\sm S,T)_K\big) \to \cdots \leqno (**)
\]
If (a)' holds, then the right hand terms in $(**)$ vanish for $i\geq 1$ in the limit, and
$(\ast)$ shows the equivalence between (b)' and (c)''.

Now assume that (b) and (c) hold. As above, (b) implies the vanishing of the middle term for $i=2$ in the limit of the sequence $(**)$. Then condition (c)' shows  (a)'. This shows that any two of the conditions (a)--(c) imply the third.

\medskip
Finally, assume that (a)--(c) hold. Tensoring, for $K\subset k_S^T(p)$, the exact sequence (see \cite{NSW}, 10.3.12)
\[
0 \to E_{K,T} \to E_{K,S\cup T} \to \bigoplus_{\p\in S(K)} (K_\p^\times/U_\p) \to \Cl_T(K) \to \Cl_{S\cup T}(K) \to 0
\]
by (the flat $\Z$-algebra) $\Z_p$, we obtain an exact sequence of finitely generated, hence compact $\Z_p$-modules. Passing to the projective limit over all finite subextensions $K$ of $k_S^T(p)|k$ and using $\varprojlim \Cl_T(K) \otimes \Z_p =0$, we obtain the exact sequence
\[
0 \to \!\!\! \varprojlim_{K\subset k_S^T(p)}\!\!\! E_{K,T} \otimes \Z_p \to \!\!\!\varprojlim_{K\subset k_S^T(p)}\!\! \! E_{K,S\cup T} \otimes \Z_p \to \!\! \varprojlim_{K\subset k_S^T(p)} \bigoplus_{\p\in S(K)} (K_\p^\times /U_\p) \otimes \Z_p \to 0.
\]
Condition (c) and local class field theory imply the vanishing of the right hand limit.
Therefore (b) implies the vanishing of the projective limit in the middle.
\end{proof}

If $G_S^T(k)(p)\neq 1$ and condition (a) of Proposition~\ref{thmb} holds, then the failure in condition (c) can only come from primes dividing~$p$. This follows from

\begin{proposition} \label{fulllocal} Let $S$ and $T$ be disjoint, finite sets of places of the global field~$k$ and let $p\neq \text{char}(k)$ be a prime number. Suppose $p\neq 2$ or $k$ totally imaginary if $k$ is a number field. If $k$ is a function field, suppose that
\[
\null_p \Cl(k)\neq 0 \ \text{ or }\ \delta=0\ \text{ or }\ T\neq\varnothing\ \text{ or }\ \# S\geq 2.
\]
If $(X\sm S,T)$ has the $K(\pi,1)$-property for~$p$ and $G_S^T(k)(p)\neq 1$, then every prime $\p\in S$ with $\zeta_p\in k_\p$ has an infinite inertia group in $G_S^T(k)(p)$. Furthermore,
\[
k_S^T(p)_\p=k_\p(p)
\]
for all  $\p\in S_\min \sm S_p$.
\end{proposition}

\begin{example} Let $\F$ be a finite field of order $\# \F \equiv 1 \bmod p$. We put $k=\F(t)$, $T=\varnothing$, and let $S$ consist of the infinite place, i.e.\ the place associated to the degree valuation. Then ${\mathbb A}^1_\F={\mathbb P}^1_{\F} \sm \{\infty\}$ has the $K(\pi,1)$-property for~$p$ and $k_{\{\infty \}}(p)|k$ is the (unramified) cyclotomic $\Z_p$-extension. Hence the assumption made for function fields in Proposition~\ref{fulllocal} is necessary.
\end{example}

\begin{proof}[Proof of Proposition~\ref{fulllocal}] Without restriction, we may assume that  $S=S_\min\neq \varnothing$. Suppose that some $\p\in S$ with $\zeta_p\in k_\p$ would be unramified in the extension $k_S^T(p)|k$. Then we have
$k_{S'}^T(p)=k_S^T(p)$ for $S'=S\sm \{\p\}$. In particular, we obtain an isomorphism
\[
H^1_\et (X\sm S', T, \F_p) \stackrel{\sim}{\lang} H^1_\et (X\sm S, T, \F_p).
\]
As before, we omit the coefficients $\F_p$ from the notation.  Using $H^3_\et(X\sm S,T)=0$, the excision sequence implies the commutative and exact diagram
\[
\xymatrix@C=.45cm{&H^2(G_{S'}^T(k)(p))\ar[r]^\sim\ar@{^{(}->}[d]&H^2(G_S^T(k)(p))\ar[d]^\wr\\
H^2_\p(X,T)\ar@{^{(}->}[r]&H^2_\et(X\sm S',T)\ar[r]^\alpha&H^2_\et(X\sm S,T)\ar[r]&H^3_\p(X,T)\ar@{->>}[r]&H^3_\et(X\sm S',T).}
\]
We see that $\alpha$ is a split surjection and that $\F_p\cong H^3_\p(X,T)\stackrel{\sim}{\to} H^3_\et(X\sm S',T)$. Theorem~\ref{globcoh} implies $S'=\varnothing$, hence $S=\{\p\}$, and $\delta=1$. As the same argument applies to every finite subextension $K$ of $k_S^T(p)|k$, $\p$ is inert in the extension $k_S^T(p)=k_\varnothing^T(p)$. Therefore the natural homomorphism
\[
\Gal(k_\p^\nr(p)|k_\p) \lang G_\varnothing^T(k)(p)
\]
is surjective, in particular, the group $G_S^T(k)(p)=G_\varnothing^T(k)(p)$ is procyclic. In the function field case we therefore have $\null_p \Cl(k)= 0$ or $T\neq \varnothing$. Because of $\# S=1$ and $\delta=1$, our assumptions imply  $T\neq \varnothing$ if $k$ is a function field. Therefore, by class field theory, the group $G_\varnothing^T(k)(p)$
is finite. As this group is nontrivial by assumption, it has infinite cohomological dimension, a contradiction to the $K(\pi,1)$-property. Hence every $\p\in S$ with $\zeta_p\in k_\p$  ramifies in $k_S^T(p)$. The same argument applies to every finite extension of $k$ inside $k_S^T(p)$, and so the inertia groups must be infinite.
For $\p\notin  S_p$ this implies $k_S(p)_\p=k_\p(p)$, as can be easily deduced from the explicitly known structure of the group $\Gal(k_\p(p)|k_\p)$  (cf.\ \cite{NSW}, 7.5.2).
\end{proof}

\begin{theorem} \label{dualmod} Let $S\neq \varnothing$ and $T$ be disjoint, finite sets of places of the global field~$k$ and let $p\neq \text{char}(k)$ be a prime number. Suppose $p\neq 2$ or $k$ totally imaginary if $k$ is a number field. We assume that conditions {\rm (a)--(c)} of Proposition~\ref{thmb} hold and that  $\zeta_p\in k_\p$ for all $\p\in S$.

\smallskip
Then $G_S^T(k)(p)$ is a pro-$p$ duality group of dimension~$2$.
\end{theorem}

\begin{proof} By condition (a), we have  $H^3(G_{S}^T(k)(p),\F_p) \stackrel{\sim}{\to} H^3_\et (X\sm S, T, \F_p)=0$, hence $\cd\, G_{S}^T(k)(p)\leq 2$. On the other hand, by (c), the group $G_{S}^T(k)(p)$ contains the full local group $\Gal(k_\p(p)|k_\p)$ as a subgroup for each $\p\in S$.  Because of $\zeta_p\in k_\p$ for $\p\in S$, these local groups have cohomological dimension~$2$, and so also $\cd\,G_{S}^T(k)(p)=2$.

By \cite{NSW}, Theorem~3.4.6, in order to show that $G_S^T(k)(p)$ is a duality group, we have to prove the vanishing of the terms
\[
D_i\big(G_S^T(k)(p)\big): = \varinjlim_{\substack{U\subset G_S^T(k)(p)\\ \cor^\vee}} H^i(U, \F_p)^\vee,\quad i=0,1,
\]
where $U$ runs through the open normal subgroups of $G_S^T(k)(p)$ and the transition maps are the duals of the corestriction homomorphisms. The vanishing of $D_0$ is obvious as $G_S^T(k)(p)$ is infinite. We therefore have to show that
\[
\varinjlim_{K\subset k_S^T(p)} H^1\big((X\sm S,T)_K,\F_p\big)^\vee=0.
\]
We assume for the moment that $k$ is a number field or  $T\neq \varnothing$. In this case, for every finite extension  $K$ of $k$ inside $k_S^T(p)$,  the group $\Cl_T(K)$ is finite and the principal ideal theorem implies
that
\[
\varinjlim_{K\subset k_S^T(p)} H^1_\et\big((X,T)_K,\F_p\big)^\vee= \varinjlim_{K\subset k_S^T(p)} \Cl_T(K)/p =0.
\]
Excision and the local duality theorem imply the exactness of the sequence
\[
\bigoplus_{\p\in S(K)} H^1_\nr(K_\p,\mu_p) \to H^1\big((X\sm S,T)_K,\F_p\big)^\vee \to H^1\big((X,T)_K,\F_p\big)^\vee\,.
\]
Since $\zeta_p\in k_\p$ and  $k_S^T(p)_\p=k_\p(p)$ for all $\p\in S$, the left hand term, vanishes in the limit over all $K\subset k_S^T(p)$. We have seen above that also the right hand term vanishes in the limit. This shows the statement.
It remains to deal with the case $T=\varnothing$ if $k$ is a function field. By Poincar\'{e} duality, we have
\[
H^1_\et\big((X\sm S)_K,\F_p\big)^\vee \cong H^2_c\big((X\sm S)_K,\mu_p\big)= H^2_\et(X_K,j_!\mu_p),
\]
where $j: (X\sm S)_K \to X_K$ denotes the inclusion. The excision sequence together with $H^2_\p(X_K,j_!\mu_p)\cong H^1(K_\p,\mu_p)$ for $\p\in S(K)$ shows the exactness of the sequence
\[
\bigoplus_{\p \in S(K)} H^1(K_\p,\mu_p) \to H^2_\et(X_K, j_! \mu_p) \to H^2_\et\big((X\sm S)_K, \mu_p\big).\leqno (\ast)
\]
Again because of $\zeta_p\in k_\p$ and $k_S(p)_\p=k_\p(p)$ for $\p\in S$, the left hand term in $(\ast)$ vanishes in the limit over all $K\subset k_S(p)$.  The Kummer sequence implies the exactness of
\[
\Cl_S(K)/p \lang H^2_\et((X\sm S)_K, \mu_p) \lang \null_p \Br((X\sm S)_K).\leqno (\ast\ast)
\]
The principal ideal theorem shows  $\varinjlim_{K\subset k_S(p)} \Cl_S(K)/p=0$ and the Hasse principle for the Brauer group provides an injection
\[
\null_p \Br\big((X\sm S)_K\big) \hookrightarrow \bigoplus_{\p\in S(K)} \null_p\Br(K_\p).
\]
Since $k_S(p)$ realizes the maximal unramified $p$-extension of  $k_\p$ for $\p\in S$, the right hand, and hence also the middle term of  $(\ast\ast)$ vanishes in the limit over  $K$. Hence also the middle term of  $(\ast)$ vanishes in the limit. This finishes the proof.
\end{proof}

\begin{remark} For $S\neq \varnothing$, $S$-$T$-id\`{e}le class group of $k_S^T$ is defined by
\[
C_S^T:=\varinjlim_{K\subset k_S^T} \coker \big(E_{K,S\cup T} \lang \prod_{\p\in S(K)} K_\p^\times\big)\,.
\]
The pair $(G_S^T(k), C_S^T)$ is a class formation.  For a subfield $K\subset k_S^T$, one defines $C_S^T(K):=(C_S^T)^{\Gal(k_S^T|K)}$.
Under the assumptions of Theorem~\ref{dualmod}, the dualizing module of the duality group $G_S^T(k)(p)$ is isomorphic to $\text{\rm tor}_{p}\big(C_S^T(k_S^T(p))\big)$.
\end{remark}

\bigskip

\noindent {NWF I - Mathematik, Universit\"{a}t Regensburg, D-93040
Regensburg, Deutschland. Email-address: alexander.schmidt@mathematik.uni-regensburg.de}

\end{document}